\documentclass[journal, twoside]{IEEEtran}

\usepackage{import}

\usepackage[utf8]{luainputenc}

\usepackage[T1]{fontenc}

\usepackage{graphicx}
\usepackage[dvipsnames]{xcolor}

\usepackage[english]{babel}
\addto\captionsenglish{\renewcommand{\figurename}{Fig.}}
\addto\captionsenglish{\renewcommand{\tablename}{Table}}
\usepackage{csquotes}

\usepackage{newtxtext}
\usepackage{amsthm}
\usepackage[slantedGreek]{newtxmath}
\usepackage[OMLmathsfit]{isomath}
\DeclareMathAlphabet{\mathbfsf}{\encodingdefault}{\sfdefault}{bx}{n}
\usepackage{bm}
\usepackage{envmath}
\usepackage{mathtools}
\usepackage{commath}
\usepackage{siunitx}
\usepackage{nicefrac}

\usepackage[caption=false,font=footnotesize]{subfig}

\usepackage{booktabs}
\usepackage{footmisc}  

\usepackage{url}

\theoremstyle{definition}

\theoremstyle{plain}

\theoremstyle{remark}

\usepackage{lineno}
\modulolinenumbers[5]
\usepackage{todonotes}
\usepackage{umoline}

\usepackage{pgfplots}
\usepackage{pgfplotstable}
\pgfplotsset{compat=newest}
\pgfplotsset{plot coordinates/math parser=false}
\newlength\figureheight
\newlength\figurewidth
\pgfplotsset{every axis plot/.append style={line width=1.5pt},
    legend style={font=\footnotesize, 
        text height=1.0ex,
        draw=black,
        fill=white,
        legend cell align=left}}

\usepackage[hidelinks]{hyperref} 
\usepackage[english]{cleveref}

\Crefname{defn}{definition}{definitions}
\Crefname{defn}{Definition}{Definitions}

\Crefname{asm}{assumption}{assumptions}
\Crefname{asm}{Assumption}{Assumptions}

\crefname{lem}{lemma}{lemmas} 
\Crefname{lem}{Lemma}{Lemmas}

\crefname{prop}{proposition}{propositions} 
\Crefname{prop}{Proposition}{Propositions}

\crefname{thm}{theorem}{theorms} 
\Crefname{thm}{Theorem}{Theorms}

\crefname{cor}{corollary}{corollaries}
\Crefname{cor}{Corollary}{Corollaries}
\newcounter{subequation}
\newlength\mtabskip\mtabskip=-1.25cm

\def\mtabLong{long}

\newcommand{\mr}{\mathrm}

\newcommand{\veg}[1]{\bm{#1}}     
\renewcommand{\vec}[1]{\mathsfbfit{#1}} 
\newcommand{\vecop}[1]{\bm{\mathcal{#1}}} 

\newcommand{\matel}[1]{\begin{bmatrix} #1 \end{bmatrix}}

\newcommand{\matI}{\mathbfsf{I}}

\newcommand{\n}{\hat{\bm{n}}}

\newcommand{\dd}{\mathrm{d}}  

\newcommand{\im}{\mathrm{i}}  

\newcommand{\e}{\mathrm{e}}

\newcommand{\T}{\mr{T}}

\newcommand\restr[2]{{
        \left.\kern-\nulldelimiterspace 
        #1 
        \vphantom{|} 
        \right|_{#2} 
}}

\newcommand\rst[3]{{
        \left.\kern-\nulldelimiterspace 
        #1 
        \vphantom{|} 
        \right|_{#2}^{#3} 
}}

\newcommand{\TA}{\vec T_{\!\!\mr{A}}}

\newcommand{\nT}{\,{}^{\vec n}\!\mskip2mu\vec T} 
\newcommand{\nTA}{\,{}^{\vec n}\!\mskip2mu\vec T_{\!\!\mr{A}}}
\newcommand{\nK}{\,{}^{\vec n}\!\!\mskip2mu \vec K} 
\newcommand{\nR}{\,{}^{\vec n}\!\!\mskip2mu \vec R} 

\newcommand{\nev}{\,{}^{\vec n}\!\vec e}
\newcommand{\nhv}{\,{}^{\vec n}\!\vec h} 
\usepackage{acro}

\DeclareAcronym{DG}
{
    short = DG ,
    long = discontinuous Galerkin
}

\DeclareAcronym{ACA}
{
    short = ACA ,
    long = adaptive cross approximation
}

\DeclareAcronym{EFIE}
{
    short =  EFIE ,
    long = electric field integral equation
}

\DeclareAcronym{MFIE}
{
    short =  MFIE ,
    long = magnetic field integral equation
}

\DeclareAcronym{CFIE}
{
    short =  CFIE ,
    long = combined field integral equation
}

\DeclareAcronym{EFIO}
{
    short =  EFIO ,
    long = electric field integral operator
}

\DeclareAcronym{MFIO}
{
    short =  MFIO ,
    long = magnetic field integral operator
}

\DeclareAcronym{CFIO}
{
    short =  CFIO ,
    long = combined field integral operator
}

\DeclareAcronym{MUIE}
{
    short =  MUIE ,
    long = Müller integral equation
}

\DeclareAcronym{PMCHWT}
{
    short =  PMCHWT ,
    long = Poggio-Miller-Chang-Harrington-Wu-Tsai integral equation
}

\DeclareAcronym{SPD}
{
    short =  SPD ,
    long = {symmetric, positive definite}
}

\DeclareAcronym{SPSD}
{
    short =  SPD ,
    long = {symmetric, positive semi-definite}
}

\DeclareAcronym{PEC}
{
    short =  PEC ,
    long = perfectly electrically conducting
}

\DeclareAcronym{RWG}
{
    short = RWG ,
    long = Rao-Wilton-Glisson
} 

\DeclareAcronym{BC}
{
    short = BC ,
    long = Buffa-Christiansen
}

\DeclareAcronym{SVD}
{
    short = SVD ,
    long = singular value decomposition
}

\DeclareAcronym{CG}
{
    short = CG ,
    long = conjugate gradient
} 

\DeclareAcronym{PCG}
{
    short = PCG ,
    long = preconditioned conjugate gradient
} 

\DeclareAcronym{CGS}
{
    short = CGS ,
    long = conjugate gradient squared
}

\DeclareAcronym{CMP}
{
    short = CMP ,
    long = Calderón multiplicative preconditioner
} 

\DeclareAcronym{RFCMP}
{
    short = RF-CMP ,
    long = refinement-free Calderón multiplicative preconditioner
} 

\DeclareAcronym{HPD}
{
    short = HPD ,
    long = {Hermitian, positive definite}
} 

\DeclareAcronym{RHS}
{
    short = RHS ,
    long = right-hand side
}

\DeclareAcronym{LSE}
{
    short = LSE ,
    long = linear system of equations
}

\DeclareAcronym{AMG}
{
    short = AMG ,
    long = algebraic multigrid
}

\DeclareAcronym{PW}
{
    short = PW ,
    long = plane wave
}

\DeclareAcronym{GMRES}
{
    short = GMRES ,
    long = generalized minimum residual
}

\DeclareAcronym{IDR}
{
    short = IDR ,
    long = induced dimension reduction
}

\DeclareAcronym{BICGstab}
{
    short = BiCGstab ,
    long = stabilized bi-conjugate gradient
}

\DeclareAcronym{FF}
{
    short = FF ,
    long = far field
}

\DeclareAcronym{NF}
{
    short = NF ,
    long = near field
}

\DeclareAcronym{TE}
{
    short = TE,
    long = transverse electric
}

\DeclareAcronym{TM}
{
    short = TM,
    long = transverse magnetic
}

\DeclareAcronym{NURBS}
{
    short = NURBS,
    long = non-uniform rational B-splines
}

\DeclareAcronym{FEM}
{
    short = FEM,
    long = finite element method
}

\DeclareAcronym{PEEC}
{
    short = PEEC,
    long = partial element equivalent circuit
}

\DeclareAcronym{GWP}
{
    short = GWP,
    long = Graglia Wilton Peterson
}

\DeclareAcronym{MoM}
{
    short = MoM,
    long = method of moments
}

\DeclareAcronym{GPDIM}
{
    short = GP-DIM,
    long = general-purpose density interpolation method
} 

\newcolumntype {n}{c}
\newcolumntype {N}{>{\small}c}
\newcolumntype {L}{>{\small}l}
\newcolumntype {F}{>{\footnotesize}c}
\newcolumntype {v}[1]{>{\raggedright \hspace {0pt}} p {#1}}
\newcolumntype {V}[1]{>{\small \raggedright \hspace {0pt}} p {#1}}
\newcolumntype{d}[1]{>{\DC@{.}{.}{#1}}c<{\DC@end}}

\newcolumntype{R}[1]{%
    >{\begin{turn}{90}\begin{minipage}{#1}\small\raggedright\hspace{0pt}}l%
            <{\end{minipage}\end{turn}}%
}

\usepackage[style=ieee, backend=biber, isbn=false,maxbibnames=99]{biblatex}
  
\usepackage{multirow}

\graphicspath{{figures/}} 
\usepackage{atbegshi}
\usepackage{lipsum}

\begin{document}

    \title{High-Order-Accurate Continuity Enforcing \\Nyström Discretization of 3D Maxwell \\Combined Field Integral Equations}

	%

    \author{Bernd Hofmann,~\IEEEmembership{Member,~IEEE,}
            Reza Molavi,~\IEEEmembership{Member,~IEEE,}
	        and~Constantine Sideris,~\IEEEmembership{Senior Member,~IEEE\vspace{-2mm}}
	\thanks{This work was supported by Google Quantum AI, the Air Force Office of Scientific Research (FA9550-26-1-B177, FA9550-25-1-0020), and the National Science Foundation (CCF-2047433).}%
	\thanks{B. Hofmann and C. Sideris are with the Department of Electrical Engineering, Stanford University, Stanford, CA 94305, USA (e-mail: bernd.hofmann@stanford.edu).}
	\thanks{R. Molavi is with Google Quantum AI, Santa Barbara, CA 93111, USA.}%
    \thanks{Digital Object Identifier}%
	}

	%
	%

	\markboth{}
	{Hofmann \MakeLowercase{\textit{et al.}}: High-Order-Accurate Continuity Enforcing Nyström Discretization of 3D Maxwell CFIEs}
	%



	\maketitle

    \begin{abstract}
In Nyström-collocation discretizations of the \ac{EFIE}, the surface divergence acts on surface densities that may be discontinuous across patch boundaries, which degrades accuracy and convergence. 
We show that this not only affects the \ac{EFIE} but every formulation in which the operator occurs, either in the equation itself or in the scattered field computation, and propose a high-order-accurate continuity-enforcing scheme for smooth surfaces as a remedy for the direct and indirect \acp{EFIE}, \acp{MFIE}, and regularized \acp{CFIE} alike.
The scheme comprises two ingredients: 
i) We show how to discretize the equations via a Chebyshev-based Nyström scheme, which admits closed quadrature rules. 
ii) Since unknowns and test vectors are in terms of patch-local curvilinear bases, continuity is enforced by a change of basis: we construct sparse mapping matrices assembled solely from the curvilinear geometry description. 
In doing so, we restore the accuracy of the \ac{EFIE} such that it can be combined with the \acp{MFIE} with equal weights to form \acp{CFIE}.
Numerical studies for the scattering from canonical and realistic geometries show that all considered formulations individually and combined benefit from the continuity enforcement in terms of better conditioning, reduced iterations of an iterative solver, and several more digits of accuracy in the scattered fields, despite reducing the total number of unknowns.
    \end{abstract}
    \acresetall


    %
    \IEEEpeerreviewmaketitle

    \section{Introduction}

\IEEEPARstart{T}{he} direct and indirect \acp{EFIE} and \acp{MFIE} are well-established formulations to numerically solve electromagnetic radiation and scattering problems involving \ac{PEC} objects, where the direct formulations determine the physically induced surface current density, and the indirect formulations determine non-physical surface densities and relate them to the physical scattered fields~\cite{hsiaoMathematicalFoundationsError1997,brunoElectromagneticIntegralEquations2009}.
However, for closed structures, these equations are known to possess a non-trivial nullspace at the resonance frequencies of the \ac{PEC} cavity.
The corresponding non-uniqueness of the solution results in severely deteriorated accuracy and convergence of numerical solution strategies.
A classic remedy is to form \acp{CFIE} by combining both \acp{EFIE} and \acp{MFIE}~\cite{mautzHFieldEFieldCombinedField1978,brunoElectromagneticIntegralEquations2009}.

Most commonly, these integral equations are discretized using low-order \ac{MoM} schemes (also called Petrov-Galerkin schemes) by triangulating the objects' surfaces and defining low-order basis functions on the triangulation. 
These basis functions then serve as trial functions, expanding the unknown surface densities and testing functions.
Consequently, costly 4D integrals have to be evaluated to obtain the matrix entries of the \ac{LSE} to be solved to obtain the expansion coefficients.
Higher-order versions of the \ac{MoM} discretization increase the cost of these 4D integrals even further~\cite{gragliaHigherOrderInterpolatory1997,jorgensenHigherOrderHierarchical2004,valdesHighorderDivQuasi2011}.
Moreover, an \ac{MFIE} which matches the accuracy of the \ac{EFIE} and does not deteriorate the accuracy of the \ac{CFIE} has so far only been obtained by employing a dual basis, which increases the computational effort further~\cite{coolsAccurateConformingMixed2011,begheinSpaceTimeMixedGalerkin2013,valdesHighorderDivQuasi2011,hofmannKleinheubach2024,hofmannExplicitHigherOrderDual2026}.

In contrast, we are considering an inherently high-order-accurate Nyström-collocation discretization of the integral equations, in which obtaining the matrix entries requires, in essence, only an integral kernel evaluation. 
Specifically, we are considering the approach proposed in~\cite{huChebyshevBasedHighOrderAccurateIntegral2021,garzaBoundaryIntegralEquation2020,garzaFastInverseDesign2023}, which uses Chebychev-polynomial based interpolation in combination with a singularity cancellation approach to address occurring kernel singularities in an efficient and highly accurate manner.
Other approaches to address the kernel singularities, such as the locally corrected Nyström method~\cite{caninoNumericalSolutionHelmholtz1998a,gedneyDerivingLocallyCorrected2003,gedneyLocallyCorrectedNystrom2014} or the \ac{GPDIM}, require the solution of an \ac{LSE} or a least-squares problem per patch, which can limit accuracy and efficiency.

\begin{table*}[t]
    \centering
    \caption{Considered integral equation formulations for the scattering from \ac{PEC} objects: the integral operators $\vecop T$, $\vecop K$, and $\vecop R_{\!w}$ relate the (physical) surface current density $\veg j$ and the (non-physical) surface densities $\veg a$, $\veg b$, $\veg c$, to the excitation fields $(\veg e^\mr{ex}, \veg h^\mr{ex})$. }
    \label{formulations}
    \small
    \setlength{\tabcolsep}{4pt}
    \renewcommand{\arraystretch}{1.8}
    \begin{tabular}{@{}r r @{${}={}$} l r @{${}={}$} l@{}}
    \toprule
\multicolumn{1}{c}{formulation}     & \multicolumn{2}{c}{direct} & \multicolumn{2}{c}{indirect} \\
\cmidrule(r){1-1} \cmidrule(lr){2-3}           \cmidrule(l){4-5}
EFIE   & ${(\eta \vecop T \veg j)}_\text{tan}$     & $\veg e^\mr{ex}_\text{tan}$                
       & ${(\eta \vecop T \veg j)}_\text{tan}$     & $\veg e^\mr{ex}_\text{tan}$ \\
MFIE   & $(\pm \vecop I /2 + \n \times \vecop K)\veg j$       & $\n \times \veg h^\mr{ex}$                 
       & $(\pm \vecop I /2 - \n \times \vecop K)\veg a$       & $\n \times \veg e^\mr{ex}$ \\
CFIE   & $(\pm \vecop I /2 + \n \times \vecop K) \veg j +  {( \vecop T \veg j)}_\text{tan}$          & $\n \times \veg h^\mr{ex} + \eta^{-1}\veg e^\mr{ex}_\text{tan}$ 
       & $(\pm \vecop I /2 - \n \times \vecop K) \veg b +  k\n \times \vecop T (\n \times \veg b)$    & $\n \times \veg e^\mr{ex}$  \\
CFIE-R & $(\pm \vecop I /2 + \n \times \vecop K) \veg j -  k\n \times \vecop R_{\!w}(\n \times \vecop T \veg j)$ & $\n \times \veg h^\mr{ex} -  \eta^{-1} k\n \times \vecop R_{\!w}(\n \times \veg e^\mr{ex})$ 
       & $(\pm \vecop I /2 - \n \times \vecop K) \veg c -  k\n \times \vecop T(\n \times \vecop R_{\!w} \veg c)$ & $\n \times \veg e^\mr{ex}$  \\
    \bottomrule
    \end{tabular}
\end{table*}

Converse to the \ac{MoM} discretizations, Nyström approaches have been found to yield significantly higher accuracies with the \ac{MFIE} than with the \ac{EFIE}~\cite{petersonObservedBaselineConvergence2008,coolsAccurateConformingMixed2011,gedneyDerivingLocallyCorrected2003}.
The root cause of this discrepancy is the divergence-conformity of the surface densities: the \ac{EFIE} operator contains a surface divergence, so a density exhibiting a jump across a patch boundary produces spurious line charges there.
Divergence-conformity is a necessary and sufficient condition for the fields to be of finite energy~\cite[p.~49]{cessenatMathematicalMethodsElectromagetism1996} and is thus baked into the equations as an assumption, while the equations do not, by themselves, enforce it.
In \ac{MoM} discretizations, the assumption is honored by construction, through divergence-conforming basis functions ensuring normal continuity of the surface densities~\cite{petersonMappedVectorBasis2006}.
Nyström-collocation schemes impose no such constraint.
Notably, the \ac{EFIE} operator enters not only the \acp{EFIE} but also the representation formulae from which the scattered fields are computed once a surface density has been determined.
The accuracy deficit is therefore not confined to the formulation in which the operator appears explicitly.
As an alternative, additional line integrals can be introduced to account for the discontinuities, an approach known as discontinuous Galerkin in the \ac{MoM} 
context~\cite{pengDiscontinuousGalerkinSurface2013,martinDiscontinuousGalerkinIntegral2023}; applied to Nyström discretizations, however, it has been found to increase the number of iterations, with no clear preconditioning strategy available~\cite{gedneyDerivingLocallyCorrected2003,youngHighOrderNystromImplementation2012,huHighOrderNystromBasedScheme2024}.
A different route was taken in~\cite{hendijaniConstrainedLocallyCorrected2015,pfeifferNumericalCharacterizationDivergenceConforming2017}, where a relation between the Nyström unknowns and high-order continuity-enforcing basis functions was established, inheriting the properties of the high-order \ac{MoM} at the cost of a singular value decomposition of the constraint matrices.

In this work, we propose a scheme that explicitly enforces continuity of the surface densities in the Nyström discretization by using a closed quadrature rule and constructing sparse continuity-enforcing matrices.
To this end, we i) detail how the Nyström discretization of~\cite{huChebyshevBasedHighOrderAccurateIntegral2021,garzaBoundaryIntegralEquation2020,garzaFastInverseDesign2023} can be generalized to handle also the direct and indirect \acp{CFIE} as well as the regularized versions of the latter, introduced in~\cite{brunoElectromagneticIntegralEquations2009} to mitigate the increased solution times caused by the first-kind nature of the \acp{EFIE}.
ii) We generalize the approach of~\cite{huHighOrderNystromBasedScheme2024}, which considers the direct \ac{EFIE}, by constructing the continuity-enforcing matrices via a change of basis between the patch-local curvilinear coordinate systems. 
This allows us to handle all formulations and requires three distinct sparse mapping matrices, one for the unknowns and one each for covariant and contravariant test vectors.
All three are assembled solely from the curvilinear geometry description, are constructed once, and leave the application of the discretized operators untouched.
iii) By numerical studies of the scattering from canonical and realistic geometries, we show that not only the \acp{EFIE} but all considered formulations individually and combined in the \acp{CFIE} profit in terms of better conditioning, reduced iterations of an iterative solver, and several more digits of accuracy in the scattered fields, even though the total number of unknowns is reduced.
We further consider a mixed collocation scheme, inspired by \ac{MoM} discretizations employing primal and dual basis functions~\cite{coolsAccurateConformingMixed2011,begheinSpaceTimeMixedGalerkin2013}, which we find offers no notable benefit in the Nyström setting.
Note that preliminary results have been presented in~\cite{hofmannCFIE2026}.

This article is organized as follows: 
Section~II introduces background material about the fundamental integral equation formulation and fixes the notation. 
The generalized Nyström-collocation discretization of the \acp{EFIE}, \acp{MFIE}, and regularized \acp{CFIE} is introduced in Section~III, and in 
Section~IV, the continuity enforcement for the considered discretizations is derived.
The numerical studies are presented in Section~V, followed by the conclusion.

    \section{Background: Electric, Magnetic, and Combined Field Integral Equations}

Let a time-harmonic field $(\veg e^\mr{ex}, \veg h^\mr{ex})$ with an assumed but throughout this article suppressed time dependency of $\e^{-\im \omega t}$ excite a \ac{PEC} object with surface $\Gamma \in \mathbb{R}^3$ equipped with a unit normal vector field $\n$ embedded in a homogeneous background medium with permittivity $\varepsilon$ and permeability $\mu$. 
The surface $\Gamma$ is assumed to be a bounded, smooth, two-dimensional manifold which can be simply or multiply connected.
In order to determine the scattered field $(\veg e^\mr{sc}, \veg h^\mr{sc})$, we consider the direct and indirect integral equation formulations given in Table~\ref{formulations}, which first relate a surface density $\veg j$, $\veg a$, $\veg b$, or $\veg c$ on $\Gamma$ to the excitation fields $(\veg e^\mr{ex}, \veg h^\mr{ex})$.
All formulations are based on the identity operator $\vecop I \veg j = \veg j$, the integral operator
\begin{multline}
\vecop T \veg j = \im k \iint_\Gamma g_k(\veg r, \veg r') \veg j(\veg r') \,\dd \veg r' \\ 
+ \im k^{-1} \nabla \!\! \iint_\Gamma g_k(\veg r, \veg r') \, \nabla_\Gamma' \cdot \veg j(\veg r') \, \dd \veg r'\,
    \label{Top}
\end{multline} 
and the integral operator
\begin{equation}
\vecop K \veg j = \iint_\Gamma \nabla g_k(\veg r, \veg r') \times \veg j(\veg r') \, \dd \veg r'\,,
    \label{Kop}
\end{equation} 
involving imaginary unit $\im^2=-1$, wave impedance $\eta = \sqrt{\mu / \varepsilon}$, wavenumber $k=\omega\sqrt{\mu \varepsilon}$, and the free-space scalar Green function $g_k(\veg r, \veg r') = \e^{\im k |\veg r - \veg r'|} \big/ ({4\uppi |\veg r - \veg r'|})$.
The integral in $\n \times \vecop K \veg j$ is to be interpreted in a principal value sense (see, e.g., \cite{hsiaoMathematicalFoundationsError1997,buffagalerkin2003} for a more subtle treatment of the traces).
The direct \ac{EFIE} and \ac{MFIE}~\cite{maueZurFormulierungAllgemeinen1949} can be derived from equivalence principles and determine the physically induced surface current density $\veg j$.
The indirect \ac{EFIE} and \ac{MFIE} can be derived from an ansatz for the scattered fields~\cite{hsiaoMathematicalFoundationsError1997,coltonIntegralEquationMethods1983} and determine non-physical surface densities\footnote{Note that indirect and direct \ac{EFIE} are identical, where sometimes the indirect \ac{EFIE} is defined with a different constant prefactor. Spectrally, indirect and direct \ac{MFIE} behave essentially the same~\cite{hsiaoMathematicalFoundationsError1997,brunoElectromagneticIntegralEquations2009}.}.  
The occurring $\pm \vecop I$ correspond to exterior and interior scattering problems, i.e., the excitation fields are outside (+) or inside (-) a closed $\Gamma$.

As the indirect and direct \ac{EFIE} and \ac{MFIE} are well-known to break down at the (same) cavity resonance frequencies of $\Gamma$, \acp{CFIE} have been established as a means to overcome that breakdown.
The direct \ac{CFIE}~\cite{mautzHFieldEFieldCombinedField1978} is formed by a linear combination of the \ac{EFIE} and the \ac{MFIE}. 
The indirect \ac{CFIE} is formed by a suitable ansatz with unique solvability shown in~\cite{brunoElectromagneticIntegralEquations2009}.
Notably, only the incident electric field $\veg e^\mr{ex}$ is involved, in contrast to the direct formulation using both $\veg e^\mr{ex}$ and $\veg h^\mr{ex}$. 
Note that, in our definitions for the \acp{CFIE} in Table~\ref{formulations}, we have weighted \ac{EFIE} and \ac{MFIE} contributions equally, as we aim to discretize them in a way that provides similar accuracy for the scattered fields.

As the integral operator $\vecop T$ is of the first kind, the iterative solution of the \ac{EFIE} usually results in significantly more iterations than for the \acp{MFIE}. 
To mitigate this effect for the \acp{CFIE} a regularization approach has been proposed in~\cite{brunoElectromagneticIntegralEquations2009}: by applying the regularizer
\begin{equation}
\vecop R_{\!w} \veg j = \iint_\Gamma g_w(\veg r, \veg r') \veg j(\veg r') \,\dd \veg r' 
    \label{Rop}
\end{equation} 
to the operator $\vecop T$, the composition is compact such that in total the \acp{CFIE} are the sum of an identity and compact operators.
The unique solvability for $k>0$ is maintained by choosing $w=\im k / 2$ for the regularizer.
Notably, $\vecop R_{\!w}$ acts in the direct formulation as a left regularizer, whereas in the indirect formulation, $\vecop R_{\!w}$ acts as a right regularizer.
This could make a difference as it potentially smooths the density $\veg c$ before the surface divergence operator occurring in $\vecop T$ is applied to it. 
This has been leveraged on in~\cite{siderisHighOrderaccurateSolution2025a} for the 2D scalar Helmholtz equation, and its influence on the continuity enforcement will be investigated in the numerical results.

\begin{table}[t]
    \centering
    \caption{Computation of the scattered fields $(\veg e^\mr{sc}, \veg h^\mr{sc})$ from the surface densities.}
    \label{potentials}
    \small
    \setlength{\tabcolsep}{7pt}
    \renewcommand{\arraystretch}{1.5}
    \begin{tabular}{@{}clll@{}}
    \toprule
density           & \multicolumn{1}{c}{$\veg e^\mr{sc}$}               & \multicolumn{1}{c}{$\veg h^\mr{sc}$} \\
\cmidrule(r){1-1}   \cmidrule(lr){2-2}                                   \cmidrule(l){3-3}
$\veg j$ & $\eta \vecop T \veg j$                                           & $\vecop K \veg j$ \\
$\veg a$ & $-\vecop K \veg a$                                               & $\eta^{-1}\vecop T \veg a$ \\
$\veg b$ & $-\vecop K \veg b - k\vecop T(\n \times \veg b)$                 & $\eta^{-1}[\vecop T \veg b - k \vecop K (\n \times \veg b)]$ \\
$\veg c$ & $-\vecop K \veg c - k\vecop T(\n \times \vecop R_{\!w} \veg c)$  & $\eta^{-1}[\vecop T \veg c - k \vecop K(\n \times \vecop R_{\!w} \veg c)]$ \\
    \bottomrule
    \end{tabular}
\end{table}
Once one of the surface densities $\veg a$, $\veg b$, $\veg c$, or $\veg j$ has been determined, the radiated or scattered fields can be obtained by applying the operators $\vecop T$, $\vecop K$, and $\vecop R_{\!w}$ as given in Table~\ref{potentials}.

\begin{table*}[t]
    \centering
    \caption{Interpolation between quadrature points with $\alpha_m = 1$ for $m=0$ and $\alpha_m = 2$ otherwise. Moreover, $\beta_m = 1/2$ for $m=0; m = N-1$ and $\beta_m = 1$ otherwise.}
    \label{interpol}
    \small
    \setlength{\tabcolsep}{5pt}
    \renewcommand{\arraystretch}{2.5}
    \begin{tabular}{@{}rllll@{}}
    \toprule
\multicolumn{1}{c}{quadrature}       & \multicolumn{1}{c}{interpolation}  & \multicolumn{1}{c}{coefficients} & \multicolumn{1}{c}{quadrature points}\\
\cmidrule(r){1-1}  \cmidrule(lr){2-2}                    \cmidrule(l){3-3} \cmidrule(l){4-4}
Fejer 1\textsuperscript{st}        & $\displaystyle j_{u/v}(u,v) = \sum\limits_{m=0}^{N-1}\sum\limits_{n=0}^{N-1} \gamma_{mn} T_m(u)T_n(v)$     & $\displaystyle\gamma_{mn} = \dfrac{\alpha_m \alpha_n}{N_u N_v} \sum\limits_{q=0}^{N-1}\sum\limits_{p=0}^{N-1} T_m(z_q) T_n(z_p) j_{u/v}(z_q,z_p)$         & $z_q$ zeros of the $T_q$ \\
Clenshaw-Curtis  & $\displaystyle j_{u/v}(u,v) = \sum\limits_{m=0}^{N-1}\sum\limits_{n=0}^{N-1} \gamma_{mn} T_m(u)T_n(v)$     & $\displaystyle\gamma_{mn} = \dfrac{\alpha_m \alpha_n}{N_u N_v} \sum\limits_{q=0}^{N-1}\sum\limits_{p=0}^{N-1} \beta_q \beta_p T_m(e_q) T_n(e_p) j_{u/v}(e_q,e_p)$ & $e_q$ extrema of the $T_q$\\
    \bottomrule
    \end{tabular}
\end{table*}

    \section{Generalized Nyström-Collocation Discretization} \label{secDisc}

To find approximate solutions to the equations in Table~\ref{formulations} and compute the scattered fields according to Table~\ref{potentials}, we fundamentally follow the Chebychev-Nyström-collocation discretization scheme proposed in~\cite{garzaBoundaryIntegralEquation2020,huChebyshevBasedHighOrderAccurateIntegral2021,huHighOrderNystromBasedScheme2024}, but generalize it to the regularized \acp{CFIE}.

    \subsection{Fundamental Approach and Singularity Treatment}

The surface $\Gamma$ of the scatterer is assumed to be described by a union of non-overlapping curvilinear quadrilateral patches. 
On each patch, the density of interest, e.g., $\veg j$, is represented in terms of the covariant curvilinear base vectors $\veg e_{u/v}$ as
\begin{equation}
    \veg j(u,v) = \dfrac{1}{D}\left[j_u(u,v)\veg e_u(u,v) + j_v(u, v) \veg e_v(u,v)\right]\,, 
    \label{density}
\end{equation}
with the unknown scalar functions $j_u$ and $j_v$ and $D(u,v)$ denoting the generalized Jacobi determinant of the mapping from the parametric domain $(u,v) \in {[-1,1]}^2$ to the physical domain. 
The integrals in the operators in Table~\ref{formulations} are then approximated by numerical quadratures, for example,   
\begin{multline}
    \iint_\Gamma g_k(\veg r, \veg r') \veg j(\veg r') \,\dd \veg r' \\ 
    \approx \sum_i \sum_j w_{ij} \, g_k(\veg r, \veg r_{ij}) \left[ j_{u,ij} \veg e_{u,ij} + j_{v,ij} \veg e_{v,ij} \right]
    \label{intosum}
\end{multline}
with quadrature weights $w_{ij}$ and the current density evaluated at the quadrature points $(u_i, v_j)$, i.e., $j_{u,ij}=j_u(u_i, v_j)$ and $j_{v,ij}=j_v(u_i, v_j)$ become the unknowns to be determined.

In order to handle the singularities occurring in the integral operators $\vecop T$ and $\vecop K$, both are reformulated using vector identities:
The operator $\n \times \vecop T$ is expressed as~\cite{garzaBoundaryIntegralEquation2020}
\begin{multline}
    \n \times \vecop T \veg j = \im k \n \times \iint_\Gamma g_k \veg j(\veg r') \,\dd \veg r' \\
    - \im k^{-1} D^{-1} \left[ \veg e_u \partial_v - \veg e_v \partial_u \right]  \iint_\Gamma g_k \, \nabla_\Gamma' \cdot \veg j(\veg r') \, \dd \veg r'
    \label{nTmod}
\end{multline}
and the operator $\n \times \vecop K$ as~\cite{garzaBoundaryIntegralEquation2020} 
\begin{multline}
    \n \times \vecop K \veg j = \iint_\Gamma \n \cdot \nabla g_k \,\veg j \dd \veg r' - \n \n \cdot \iint_\Gamma \n \cdot \nabla g_k \,\veg j \dd \veg r'  \\[2mm]
            - \left[ \veg e^u \n \cdot \partial_u + \veg e^v \n \cdot \partial_v \right] \iint_\Gamma g_k \veg j \,\dd\veg r'
            \label{nKmod}
\end{multline}
leaving only weakly singular integrals with kernels $g_k$ and \mbox{$\n \cdot \nabla g_k$}. 
As we also need the operator ${(\vecop{T}\veg j)}_\mathrm{tan}$, we leverage ${(\vecop{T}\veg j)}_\mathrm{tan} =  -\n \times \n \times \vecop T \veg j$ and the relation between covariant and contravariant vectors $\veg e^{u/v}$ (see, e.g.,~\cite{petersonMappedVectorBasis2006})
\begin{align}
    D^{-1} \n \times \veg e_u &= \veg e^v \\
     -D^{-1} \n \times \veg e_v &= \veg e^u
\end{align}
to obtain
\begin{multline}
    {(\vecop{T}\veg j)}_\mathrm{tan} = \im k \iint_\Gamma g_k \veg j(\veg r') \,\dd \veg r' \Big|_\mr{tan} \\
    +\im k^{-1} \left[ \veg e^v \partial_v + \veg e^u \partial_u \right]  \iint_\Gamma g_k \, \nabla_\Gamma' \cdot \veg j(\veg r') \, \dd \veg r' \,.
    \label{Tmod}
\end{multline}
Note that by $\iint_\Gamma g\veg j(\veg r') \,\dd \veg r' \big|_\mr{tan}$ we have directly removed the normal component of $\iint_\Gamma g\veg j(\veg r') \,\dd \veg r'$ as it is not needed.

To transform equations of the form~\eqref{intosum} to an \ac{LSE}, for the unknowns $j_{u/v,ij}$ stored in a vector $\vec j$, we evaluate~\eqref{intosum} at the quadrature points $\veg r_{mn}$ and take the inner product with test vectors.
To evaluate the integral for points, where the test point $\veg r_{mn}$ is close to the source point $\veg r_{ij}$,~\cite{huChebyshevBasedHighOrderAccurateIntegral2021,aslanyanFullyAutomatedAdaptive2026} took advantage of the fact that for Fejer's first quadrature rule, the quadrature points correspond to the zeros of the Chebychev polynomials $T_q$, and via the discrete orthogonality property of the $T_q$, the integrals can be performed accurately over $T_q$ via a change of variables and related to the unknowns $j_{u/v,ij}$.

To generalize this to a closed quadrature rule, we note that at its core, the approach of~\cite{huChebyshevBasedHighOrderAccurateIntegral2021} corresponds to an interpolation of the unknown functions $j_u(u,v)$ and $j_v(u,v)$ with respect to the values at the quadrature points, that is, 
\begin{equation}
    j_{u/v}(u,v) = \sum_q \sum_p j_{u/v,qp}\, \ell_q(u) \ell_p(v)
    \label{lagrange}
\end{equation}
with the unique Lagrange interpolation polynomials $\ell$.
Insertion into~\eqref{intosum} for an evaluation point $\veg r_{mn}$ close to the source patch, the integration can be performed with high accuracy with respect to $\ell(u)$ and $\ell(v)$ and related to the $j_{u/v,ij}$ via~\eqref{lagrange}.
For the specific cases of Fejer and Clenshaw-Curtis, a direct relation to the underlying Chebychev polynomials $T_q$ can be established, as shown in Table~\ref{interpol} based on orthogonality properties of the polynomials.
These representations allow for a numerically stable evaluation~\cite{pressNumericalRecipesArt2007}.

    \subsection{Testing and Handling of Derivatives}

The next step in forming an \ac{LSE} for the unknowns is to handle the derivatives in~\eqref{nTmod},~\eqref{nKmod}, and~\eqref{Tmod} and select the test vectors.
Starting with the surface divergence in~\eqref{nTmod}, we can again leverage the representation in terms of interpolation polynomials in~\eqref{lagrange}, to relate the coefficients $j_{u/v,ij}$ to their charge counterparts\footnote{Using $\nabla \cdot \veg j = \partial_u j_{u}(u,v) + \partial_v j_{v}(u,v)$ following from the definition of $\veg j$ in~\eqref{density} and properties of the covariant basis (see, e.g.,~\cite{petersonMappedVectorBasis2006}).} $\sigma_{u,ij} = \partial_u j_{u}(u_i,v_j)$ and $\sigma_{v,ij} = \partial_v j_{v}(u_i,v_j)$ as $\vec \sigma = \vec D \vec j$ with
\begin{equation}
    \vec D = \matel{ \vec D_{uv} & & \\ & \ddots & \\ & & \vec D_{uv}} \quad \text{and} \quad \vec D_{uv} = \matel{ \vec D_u & \vec D_v}
\end{equation}
where $\vec D_u$ and $\vec D_v$ denote the dense matrices resulting from taking the derivative in~\eqref{lagrange} numerically over each patch. 
Nyström discretizing the scalar operator $\iint_\Gamma g_k \sigma \dd \veg r'$ yields\footnote{The given matrix entry expressions hold for target points $\veg r_{mn}$ well separated from the source patch; entries corresponding to near interactions are replaced by the corrected weights obtained as described in the previous subsection.}
\begin{equation}
    \left[\vec S \right]_{rc} = w_{ij} \, g_k(\veg r_{mn}, \veg r_{ij} ) 
\end{equation}
with freely selectable index mappings between the Cartesian indices $ij$ and $mn$ and the linear indices $rc$ such that when testing~\eqref{nTmod} via the contravariant vectors $D \veg e^u$ and $D \veg e^v$ one obtains
\begin{equation}
    \nT \vec j = \left[\im k\vec \nTA + \im k^{-1} \tilde{\vec D} \vec S \vec D\right] \vec j
    \label{nTdisc}
\end{equation}
as the discretized counterpart to $\n \times \vecop T \veg j$ where 
\begin{equation}
    \left[\vec \nTA \right]_{rc} = \chi_{u/v} \, w_{ij} \,\veg e_{u/v,mn} \cdot  \veg e_{u/v,ij} \, g_k(\veg r_{mn}, \veg r_{ij} ) 
\end{equation}
with $\chi_{v} = -1$ and $\chi_{u} = 1$, leveraging 
\begin{align}
    D\veg e^u \cdot (\n \times \veg a) &= - \veg e_v \cdot \veg a \\
    D\veg e^v \cdot (\n \times \veg a) &= + \veg e_u \cdot \veg a \,,
    \label{ncross}
\end{align}
and $\tilde{\vec D}$ is also block diagonal with $-\vec D_u$ and $\vec D_v$ stacked on top of each other as the diagonal blocks.
In contrast, when discretizing $\vecop T \veg j$ in~\eqref{Tmod} it seems natural to test with the covariant vectors $\veg e_{u/v}$ to obtain
\begin{equation}
    \bar{\vec T} \vec j = \left[-\im k\bar{\TA} - \im k^{-1} \bar{\vec D} \vec S \vec D\right] \vec j
\end{equation}
where we denote with a bar on top the testing with covariant vectors, $\bar{\vec D}$ is block diagonal with $\vec D_u$ and $\vec D_v$ stacked on top of each other as the diagonal blocks and
\begin{equation}
    \left[\bar{\TA} \right]_{rc} = w_{ij} \,\veg e_{u/v,mn} \cdot  \veg e_{u/v,ij} \, g_k(\veg r_{mn}, \veg r_{ij} ) \,.
\end{equation}
Notably, $\nT$ and $\bar{\vec T}$ only differ by a signed permutation of the rows, as can be seen by comparing $\n \times \vecop T$ tested with $D\veg e^{u/v}$ resulting in
\begin{align}
     &-\im k \veg e_v \cdot \iint_\Gamma g_k \veg j \,\dd \veg r' - \im k^{-1} \partial_v  \iint_\Gamma g_k \, \nabla_\Gamma' \cdot \veg j \, \dd \veg r' \\
     &+\im k \veg e_u \cdot \iint_\Gamma g_k \veg j \,\dd \veg r' + \im k^{-1} \partial_u  \iint_\Gamma g_k \, \nabla_\Gamma' \cdot \veg j \, \dd \veg r'
\end{align}
and $\vecop T$ tested with $\veg e_{u/v}$ resulting in
\begin{align}
    &+\im k \veg e_u \cdot \iint_\Gamma g_k \veg j \,\dd \veg r' + \im k^{-1}\partial_u  \iint_\Gamma g_k \, \nabla_\Gamma' \cdot \veg j \, \dd \veg r' \\
    &+\im k \veg e_v \cdot \iint_\Gamma g_k \veg j \,\dd \veg r' + \im k^{-1}\partial_v  \iint_\Gamma g_k \, \nabla_\Gamma' \cdot \veg j \, \dd \veg r' \,.
\end{align}

Alternatively, $\vecop T \veg j$ in~\eqref{Tmod} can also be tested with the contravariant vectors $D\veg e^{u/v}$, which can be expressed as
\begin{equation}
    \vec T \vec j = \left[-\im k \TA - \im k^{-1} \vec G \bar{\vec D} \vec S \vec D\right] \vec j \,.
\end{equation}
This leads to a similar
\begin{equation}
    \left[\vec \TA \right]_{rc} =  w_{ij} \,\veg e^{u/v}_{mn} \cdot  \veg e_{u/v,ij} \, g_k(\veg r_{mn}, \veg r_{ij} ) 
\end{equation}
but additionally involves the contravariant components of the metric tensor
\begin{equation}
    \vec G = \matel{ \vec G_1 & & \\ & \ddots & \\ & & \vec G_P} \quad \text{with} \quad \vec G_p = \matel{ \vec I^{uu} & \vec I^{uv} \\ \vec I^{vu} & \vec I^{vv}}
\end{equation}
where the $\vec I^{uv}$ are diagonal matrices, which have inner products $\veg e^u_{ij} \cdot \veg e^v_{ij}$ for all quadrature points on each patch $p$ as their elements.
This occurs, since the divergence contribution is naturally expressed in covariant components, which the metric converts to contravariant testing.

As the operator $\n\times\vecop K$ in~\eqref{nKmod} occurs only in combination with the identity operator, we test it with contravariant vectors $D\veg e^{u/v}$ in order to obtain a true identity matrix $\matI$ as the discrete counterpart to the operator $\vecop I$.
Consequentially, we obtain
\begin{equation}
    \nK \vec j = \left[ \vec N - \vec L\right] \vec j
\end{equation}
as discrete version of $\n\times\vecop K$, where
\begin{equation}
    \left[ \vec N \right]_{rc} = w_{ij} \, \veg e^{u/v}_{mn} \cdot \veg e_{u/v,ij}\, \n_{mn} \cdot \nabla g_k(\veg r_{mn}, \veg r_{ij} ) 
\end{equation}
captures the first contribution in~\eqref{nKmod} and
\begin{equation}
     \vec L = \vec G \vec D_{xyz} \vec S_{xyz}
\end{equation}
the second contribution with 
\begin{equation}
    \left[ \vec S_{xyz} \right]_{rc} = w_{ij} \, \veg e_{u/v,ij} \, g_k(\veg r_{mn}, \veg r_{ij} )
\end{equation}
containing for each combination of source and observation point a vector entry, and $\vec D_{xyz}$ corresponds to applying the derivative matrices $\vec D_u$ and $\vec D_v$ to each of the Cartesian components of $\vec S_{xyz}$ and taking the inner product with $\n_{mn}$.

    \subsection{Operator Compositions}

The \acp{CFIE} and, in particular, their regularized versions involve compositions of the integral operators, such as
\begin{equation}
    \n \times \vecop R_{\!w}(\n \times \vecop T \veg j) \,.
    \label{RTj}
\end{equation}
Their discretization is enabled by the observation that the discretization strategy of Section~\ref{secDisc} is closed in the following sense: the operators map a tangential vector field, represented by its components with respect to the patch-local covariant basis, onto a second tangential vector field, and testing with the contravariant vectors $D\veg e^{u/v}$ recovers precisely the components of the latter in the same representation.
Compositions can therefore be discretized by matrix products of the individual discrete operators, without any intermediate change of representation.

Concretely, $\veg y = \n \times \vecop T \veg j$ is a tangential vector field, which can be expressed as
\begin{equation}
    \veg y = \dfrac{1}{D}\left[y_u(u,v)\veg e_u(u,v) + y_v(u, v) \veg e_v(u,v)\right] \,.
\end{equation}
Its components $y_u$ and $y_v$ are obtained by taking the inner products with the contravariant vectors $D \veg e^{u/v}$, which is precisely the discretization strategy underlying~\eqref{nTdisc}, i.e., $\vec y = \nT \vec j$.
Applying $\n \times \vecop R_{\!w}$ to $\veg y$ and testing again with $D\veg e^{u/v}$, the composition~\eqref{RTj} is thus discretized as
\begin{equation}
    \nR \nT \vec j
\end{equation}
with 
\begin{equation}
    \left[\nR \right]_{rc} = \chi_{u/v} \, w_{ij} \,\veg e_{u/v,mn} \cdot  \veg e_{u/v,ij} \, g_{w}(\veg r_{mn}, \veg r_{ij} ) \,,
\end{equation}
where, as in $\nTA$, the index of the covariant vector at the target point is the one complementary to the row index, and $\chi_{u/v}$ is the corresponding sign.
Since $\vecop R_{\!w}$ contains no derivative term, $\nR$ requires no counterpart to the matrices $\vec D$ and $\tilde{\vec D}$.
Its kernel $g_{w}$ is weakly singular, so the near-interaction entries are corrected by the same procedure as for $g_k$ described in Section~III-A.

The same principle applies to the remaining compositions.
For $\n \times \vecop T(\n \times \vecop R_{\!w} \veg c)$, occurring in the indirect \ac{CFIE}-R, the order of the two matrices is simply interchanged, yielding $\nT \nR \vec c$.
Note that $\vecop R_{\!w}$ then acts as a right regularizer, smoothing the density before the surface divergence contained in $\vecop T$ is applied to it, whereas in the direct formulation it acts as a left regularizer on the already differentiated density.

The term $\n \times \vecop T(\n \times \veg b)$ of the indirect \ac{CFIE} additionally requires the components of $\n \times \veg b$ in the covariant representation.
With $\n \times \veg e_u = D\veg e^v$ and $\n \times \veg e_v = -D\veg e^u$, one obtains
\begin{equation}
    \n \times \veg b = b_u \veg e^{v} - b_v \veg e^{u} \,,
\end{equation}
so that testing with $D \veg e^{u/v}$ introduces the contravariant components of the metric tensor collected in $\vec G$, and the term is discretized as $\nT \vec G \vec b$, where the signed permutation of the two components is absorbed into $\vec G$.

Finally, the excitation term $\n \times \vecop R_{\!w}(\n \times \veg e^\mr{ex})$ of the direct \ac{CFIE}-R is discretized as $\nR \nev^\mr{ex}$ with the excitation vector $\nev^\mr{ex}$ defined in the following subsection.

It should be noted that in all compositions, the intermediate vector field is kept in the redundant representation associated with the full set of quadrature points, i.e., the continuity-enforcing matrices introduced in the following section are applied only at the outer ends of the matrix products.
Numerically, we have compared it to enforcing the continuity between the operator compositions, but have observed no advantages in doing so.

    \subsection{Excitation Vectors}

The right-hand sides of the formulations in Table~\ref{formulations} are obtained by evaluating the excitation fields at the quadrature points and testing them with the same test vectors as the corresponding operators on the left-hand side.
Testing $\veg e^\mr{ex}_\mr{tan}$ with the covariant vectors $\veg e_{u/v}$ and with the contravariant vectors $D\veg e^{u/v}$ yields
\begin{equation}
    \left[\vec e^\mr{ex}\right]_r = \veg e_{u/v,mn} \cdot \veg e^\mr{ex}(\veg r_{mn})
\end{equation}
and
\begin{equation}
    \left[\bar{\vec e}^\mr{ex}\right]_r = D_{mn}\veg e^{u/v}_{mn} \cdot \veg e^\mr{ex}(\veg r_{mn})\,, 
\end{equation}
respectively.
For the rotated excitations, testing with the contravariant vectors and leveraging~\eqref{ncross} gives
\begin{equation}
    \left[\nev^\mr{ex}\right]_r = \chi_{u/v}\,\veg e_{u/v,mn} \cdot \veg e^\mr{ex}(\veg r_{mn})
\end{equation}
and
\begin{equation}
    \left[\nhv^\mr{ex}\right]_r = \chi_{u/v}\,\veg e_{u/v,mn} \cdot \veg h^\mr{ex}(\veg r_{mn})\,,
\end{equation}
where the index of the covariant vector is complementary to the row index again.

    \section{Proposed Continuity Enforcement}

The integral equations in Table~\ref{formulations} are strictly only valid if the surface densities are divergence-conforming, i.e., their surface divergence exists in a weak sense and no line charges arise along patch boundaries.
As discussed in Section~I, this is a necessary and sufficient condition following from the requirement that $\veg e$ and $\veg h$ have finite energy~\cite[p.~49]{cessenatMathematicalMethodsElectromagetism1996}, and it is not enforced by the equations themselves.
Instead of accounting for possible discontinuities by incorporating additional line integrals, we enforce the continuity of the densities explicitly, for all formulations in Table~\ref{formulations}.

As the Nyström discretization used implies the polynomial representation of the densities given in~\eqref{lagrange}, which is smooth within each patch, discontinuities can occur only at patch boundaries.
To rule them out, we choose a closed quadrature rule that also has quadrature points on the patch boundaries.
Consequently, this leads to having unknowns and test vectors that occur multiple times, corresponding to quadrature points that coincide in physical space where two or more patches share an edge or a vertex.
On initial inspection, this renders the resulting \acp{LSE} singular.
However, this duplication is the key to enforcing continuity by relating duplicate unknowns to a set of unique ones, as detailed in the following subsections.

\begin{table*}[t]
    \centering
    \caption[Discretized integral equations for the formulations in Table~\ref{formulations}.]{Discretized integral equations corresponding to the formulations in Table~\ref{formulations}, including continuity enforcement via the sparse matrices $\vec C$, $\widetilde{\vec C}$, and $\widebar{\vec C}$ with $\vec j = \vec C^\T \vec j'$, $\vec b = \vec C^\T \vec b'$, and $\vec c = \vec C^\T \vec c'$. 
    All quadratures are chosen to be a closed Clenshaw-Curtis rule. \\
    In case no continuity is enforced, an open Fejer quadrature is employed and $\vec C = \widetilde{\vec C} = \widebar{\vec C} = \matI$\,.}
    \label{formulationsDiscrete}
    \small
    \setlength{\tabcolsep}{12pt}
    \renewcommand{\arraystretch}{1.8}
    \begin{tabular}{@{}r r @{${}={}$} l r @{} l@{}}
    \toprule
\multicolumn{1}{c}{formulation}     & \multicolumn{2}{c}{direct} & \multicolumn{2}{c}{indirect} \\
\cmidrule(r){1-1} \cmidrule(lr){2-3} \cmidrule(l){4-5}
\multirow{2}{*}{EFIE}
       & $\eta \widetilde{\vec C} \vec T \vec C^\T \vec j'$ & $\widetilde{\vec C} \vec e^\mr{ex}$
       & \multirow{2}{*}{$\eta \widetilde{\vec C} \vec T \vec C^\T \vec j'$}
       & \multirow{2}{*}{${}={}\widetilde{\vec C} \vec e^\mr{ex}$} \\
       & $\eta \widebar{\vec C} \,\bar{\vec T} \vec C^\T \vec j'$ & $\widebar{\vec C} \vec e^\mr{ex}$
       & & \\
MFIE   & $\widebar{\vec C}\left[ \pm \matI /2 + \nK \right]\vec C^\T \vec j'$ & $\widebar{\vec C}\nhv^\mr{ex}$
       & $\widebar{\vec C}\left[ \pm \matI /2 - \nK \right]\vec C^\T \vec j'$ & ${}={}\widebar{\vec C}\nev^\mr{ex}$ \\
\multirow{2}{*}{CFIE}
       & $\widebar{\vec C}\left[\pm \matI /2 + \nK + \bar{\vec T} \right]\vec C^\T \vec j'$ & $\widebar{\vec C} \left[\nhv^\mr{ex} + {\eta}^{-1}\bar{\vec e}^\mr{ex}\right]$
       & \multirow{2}{*}{$\widebar{\vec C} \left[\pm \matI /2 - \nK + k\nT \vec G \right] \vec C^\T \vec b'$}
       & \multirow{2}{*}{${}={}\widebar{\vec C}\nev^\mr{ex}$} \\
       & $\left[\widebar{\vec C} (\pm\matI/2 + \nK) + \widetilde{\vec C}\vec T\right] \vec C^\mr{T}\vec{j}'$ & $\widebar{\vec C} \nhv^\mr{ex} + {\eta}^{-1}\widetilde{\vec C} \vec e^\mr{ex}$
       & & \\
CFIE-R & $\widebar{\vec C} \left[\pm \matI /2 + \nK -  k\nR \nT \right] \vec C^\T\vec j'$ & $\widebar{\vec C}\left[\nhv^\mr{ex} -  {\eta}^{-1}k\nR \nev^\mr{ex}\right]$
       & $\widebar{\vec C} \left[\pm \matI /2 - \nK -  k\nT \nR \right] \vec C^\T\vec c'$ & ${}={}\widebar{\vec C}\nev^\mr{ex}$ \\
    \bottomrule
    \end{tabular}
\end{table*}

\begin{figure}[tp]
	\centering
	\includegraphics[]{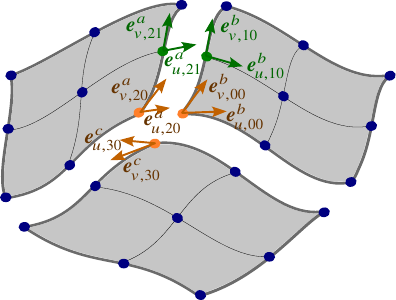}
	\captionsetup[subfloat]{labelformat=empty}
	\caption{Covariant basis vectors corresponding to 3 x 3 quadrature points on three patches sharing edges and a vertex.}
	\label{contFig}
\end{figure}

    \subsection{Surface Density}
    
Consider a setup where two patches $a$ and $b$ meet (with a smooth transition) as illustrated in Fig.~\ref{contFig} (ignoring at first the lower patch $c$): For a single point $\veg p$ on the boundary of patch $a$ with tensor product index $ij$ and its coinciding point on the boundary of patch $b$ with index $kl$, we have the density representations (exemplified for $\veg j$)
\begin{equation}
    \veg j(\veg p) = \dfrac{1}{D^{a}_{ij}} \left[j_{u,ij}^a \veg e_{u,ij}^a + j_{v,ij}^a \veg e_{v,ij}^a \right] 
\end{equation}
on patch $a$, and
\begin{equation}
    \veg j(\veg p) = \dfrac{1}{D^{b}_{kl}} \left[j_{u,kl}^b \veg e_{u,kl}^b + j_{v,kl}^b \veg e_{v,kl}^b \right] 
\end{equation}
on patch $b$; representing the same quantity in terms of different coordinate systems.
Equating and taking the inner product with the contravariant component $D^b_{kl} \veg e_{kl}^{b,u}$, we can express $j_{u,kl}^b$ in terms of $j_{u,ij}^a$ and $j_{v,ij}^a$ as
\begin{equation}
    j_{u,kl}^b = \alpha_{u,ij} \, j_{u,ij}^a + \alpha_{v,ij} \, j_{v,ij}^a
    \label{alphaa}
\end{equation}
with the coefficients
\begin{equation}
    \alpha_{u,ij} = \dfrac{D^b_{kl}}{D^a_{ij}}\, \veg e^{b,u}_{kl} \cdot \veg e^a_{u,ij}
\end{equation}
and
\begin{equation}
    \alpha_{v,ij} = \dfrac{D^b_{kl}}{D^a_{ij}}\, \veg e^{b,u}_{kl} \cdot \veg e^a_{v,ij} \,.
\end{equation}
Analogously, we can express $j_{v,kl}^b$ in terms of $j_{u,ij}^a$ and $j_{v,ij}^a$ as
\begin{equation}
    j_{v,kl}^b = \beta_{u,ij} \, j_{u,ij}^a + \beta_{v,ij} \, j_{v,ij}^a
\end{equation}
with the coefficients
\begin{equation}
    \beta_{u,ij} = \dfrac{D^b_{kl}}{D^a_{ij}}\, \veg e^{b,v}_{kl} \cdot \veg e^a_{u,ij}
\end{equation}
and
\begin{equation}
    \beta_{v,ij} = \dfrac{D^b_{kl}}{D^a_{ij}}\, \veg e^{b,v}_{kl} \cdot \veg e^a_{v,ij} \,.
    \label{betab}
\end{equation}
Hence, we can choose the two unique unknowns $j_{u,ij}^a$ and $j_{v,ij}^a$ for point $\veg p$ and obtain its corresponding values $j_{u,kl}^b$ and $j_{v,kl}^b$ via the coefficients $\alpha$ and $\beta$, which are solely determined by the different patch parametrizations.
The same principle applies when more than two patches meet at a point, as exemplified by the point where the three patches meet in Fig.~\ref{contFig}: we choose two unique unknowns $j_{u,ij}^a$ and $j_{v,ij}^a$ and compute the corresponding values on all other patches.
As~\eqref{alphaa}--\eqref{betab} express an exact change of basis, the resulting representation of $\veg j$ is independent of which patch is chosen as reference.

Identifying the coinciding points requires the connectivity of the patches, which is either directly available from the CAD description or obtained by a proximity search over the boundary quadrature points.
In addition, two adjacent patches may traverse a shared edge in opposite directions or with interchanged parametric directions, such that the index map $ij \mapsto kl$ involves a reversal and/or a transposition, which is determined once per shared edge.

Given these relations, and defining the vector of unique unknowns as $\vec j'$, it can be related to the vector with all unknowns (including all multiple occurring ones) $\vec j$ as
\begin{equation}
    \vec j = \vec C^\T \vec j'
\end{equation}
via a wide, sparse matrix 
\begingroup
\setlength{\arraycolsep}{12pt}
\renewcommand{\arraystretch}{0.35}
\begin{equation}
    \vec C = 
\begin{bmatrix} 
~\tilde{\matI} &                &        &  \\ 
               & \tilde{\matI}  &        &  \\
               &                & \ddots &  \\
               &                &        &  \tilde{\matI}~~ \\[1mm] 
      \cline{2-3} \\[-0mm]
               &               & \hspace{-6mm}\vec C_{\alpha\beta}   &     
 \end{bmatrix} \in \mathbb{R}^{N' \times N}\,,
 \label{Cmat}
\end{equation}
\endgroup
where $N' = 2P(N_u\!-\!2)(N_v\!-\!2) + 2U < N$, with $P$ the number of patches, $N_u=N_v$ the number of quadruature points in each direction, and $U$ denoting the overall number of unique patch boundary points.
The matrices $\tilde{\matI}$ denote identity matrices with all rows removed which correspond to the boundary points. 
The submatrix $\vec C_{\alpha\beta} \in \mathbb{R}^{2U \times N}$ is also sparse and contains the coefficients $\{0, 1, \alpha_i, \beta_i\}$ according to~\eqref{alphaa}-\eqref{betab}.
More precisely, per row of $\vec C_{\alpha\beta}$, there is one entry equal to 1 for the reference patch and two additional entries for each further patch touching the corresponding coinciding point.
The total number of nonzero entries is thus $O(N)$, and all of them are obtained from the description of the curvilinear geometry alone.

Note that the enforced continuity applies to the density itself, not to its surface divergence: the derivative matrices $\vec D_u$ and $\vec D_v$ in~\eqref{lagrange} remain patch-local, so $\nabla_\Gamma \cdot \veg j$ is still defined patch-wise and may jump across patch boundaries.
This is consistent with the divergence-conformity requirement, which needs a square-integrable, not a continuous, surface divergence.

    \subsection{Collocation Scheme}

In order to reduce also the number of testing vectors corresponding to the redundant set of quadrature points to (twice) the number of unique points, we could simply choose to use the two test vectors from one of the patches. 
However, to average out numerical inaccuracies (which can cause the curvilinear basis vectors of the patches to not perfectly match), we propose to follow a scheme similar to the one used for the density representation.
That is, for a contravariant test vector on patch $a$
\begin{equation}
    \veg t^a(\veg p) = D^a_{ij}  \veg e^{a,u}_{ij} \,, 
\end{equation}
we construct the corresponding contravariant test vector on patch $b$ as
\begin{equation}
    \veg t^b(\veg p) = D^b_{kl} \left[ \bar{\alpha}_{u,ij} \,\veg e^{b,u}_{kl} + \bar{\alpha}_{v,ij} \, \veg e^{b,v}_{kl} \right]\,, 
\end{equation}
with the coefficients 
\begin{equation}
    \bar{\alpha}_{u,ij} = \dfrac{D^a_{ij}}{D^b_{kl}} \,   \veg e^b_{u,kl} \cdot \veg e^{a,u}_{ij} 
\end{equation}
and 
\begin{equation}
    \bar{\alpha}_{v,ij} = \dfrac{D^a_{ij}}{D^b_{kl}}  \,  \veg e^b_{v,kl} \cdot \veg e^{a,u}_{ij} 
\end{equation}
relating the curvilinear coordinate systems of the two patches.
The actual test vector is then chosen to be the average
\begin{equation}
    \veg t(\veg p) = \dfrac{1}{n} \left[\veg t^a(\veg p) + \veg t^b(\veg p) + \veg t^c(\veg p) + \dots \right]
\end{equation}
for $n$ patches touching in the point $\veg p$.
Since $\veg t^b(\veg p)$ is by construction the expansion of $\veg t^a(\veg p)$ in the basis of patch $b$, all contributions to this average represent the same vector, and the scheme is therefore exact and independent of the reference patch; the averaging affects only the treatment of small inconsistencies in the geometry description.
The analogous scheme is applied to the test vector $\veg t^a(\veg p) = D^a_{ij}  \veg e^{a,v}_{ij}$ and corresponding coefficients $\bar{\beta}_{u,ij}$ and $\bar{\beta}_{v,ij}$. 
The mapping between all test vectors (including all vectors sharing the same quadrature point) and the unique set (with two vectors per quadrature point), can then be expressed by a sparse, wide matrix
\begingroup
\setlength{\arraycolsep}{12pt}
\renewcommand{\arraystretch}{0.35}
\begin{equation}
    \widebar{\vec C} = 
\begin{bmatrix} 
~\tilde{\matI} &                &        &  \\ 
               & \tilde{\matI}  &        &  \\
               &                & \ddots &  \\
               &                &        &  \tilde{\matI}~~ \\[1mm] 
      \cline{2-3} \\[-0mm]
               &               & \hspace{-6mm}\widebar{\vec C}_{\bar{\alpha}\bar{\beta}}   &     
 \end{bmatrix} \in \mathbb{R}^{N' \times N} \,,
 \label{Cbarmat}
\end{equation}
\endgroup
which exhibits the same structure as~\eqref{Cmat} but with a submatrix $\widebar{\vec C}_{\bar{\alpha}\bar{\beta}}  \in \mathbb{R}^{2U \times N}$ with entries $\{0, 1/n_i, \bar\alpha_i/n_i, \bar\beta_i/n_i\}$ for $n_i$ patches meeting at a point.

In case of testing with covariant vectors $\veg t^a(\veg p) = \veg e^{a}_{u,ij}$ and $\veg t^a(\veg p) = \veg e^{a}_{v,ij}$, a third kind of continuity enforcing matrix $\widetilde{\vec C}$ is needed.
It has the same structure as~\eqref{Cmat} and~\eqref{Cbarmat} but with a sparse submatrix $\widetilde{\vec C}_{\alpha\beta}$ containing the coefficients  $\{0, 1/n_i, \tilde{\alpha}_i/n_i, \tilde{\beta}_i/n_i\}$ corresponding to the coefficients in~\eqref{alphaa}-\eqref{betab} without the Jacobian determinants $D_{ij}^a$ and divided by the number of patches $n_i$ touching in a point.

An important property of the constructed matrices is
\begin{equation}
    \widebar{\vec C} \, \vec C^\T = \matI \,,
    \label{CbarCT}
\end{equation}
which follows from the fact that testing a density expanded in the covariant basis with the contravariant test vectors returns its components directly: for a coinciding point, each of the $n$ patches contributes the same value, so that the average reproduces it exactly.
As will be seen in the following subsection, this ensures that the identity operator of the \acp{MFIE} and \acp{CFIE} is mapped to a true identity matrix also after continuity enforcement, and the second-kind nature of these formulations is preserved.
Note that no analogous relation holds for $\widetilde{\vec C}$, as covariant testing of a covariantly expanded density involves the metric tensor; $\widetilde{\vec C}$ consequently appears only in combination with the \ac{EFIE} operator, which contains no identity contribution.

    \subsection{Continuity Enforced Systems}

Given the continuity enforcing matrices $\vec C$, $\widetilde{\vec C}$, and $\widebar{\vec C}$, we form the continuity enforced discretized versions of the integral equations in Table~\ref{formulations} as shown in Table~\ref{formulationsDiscrete} using the discretized operators $\vec T$, $\bar{\vec T}$, $\nT$, $\nK$, $\nR$, and $\vec G$, as well as the excitation vectors $\vec e^\mr{ex}$, $\bar{\vec e}^\mr{ex}$, $\nev^\mr{ex}$, and $\nhv^\mr{ex}$ according to their definitions in Section~\ref{secDisc}.
The pairing follows a single rule: $\vec C^\T$ is applied on the unknown side of every equation, while on the test side the matrix corresponding to the test vectors used to discretize the respective operator is applied, i.e., $\widebar{\vec C}$ for contravariant and $\widetilde{\vec C}$ for covariant testing.
The same matrix is applied to the corresponding excitation vector on the right-hand side.
All resulting systems are square of size $N' \times N'$ with $N' = 2P(N_u\!-\!2)(N_v\!-\!2) + 2U < N$.
Importantly, a closed Clenshaw-Curtis quadrature rule is chosen everywhere.
To draw comparisons to a Nyström discretization where no continuity is enforced, we also consider the case of an open Fejer quadrature rule (of the first kind) and set $\vec C = \widetilde{\vec C} = \widebar{\vec C} = \matI$ in this case.

By virtue of~\eqref{CbarCT}, the identity contributions of the \acp{MFIE} and \acp{CFIE} reduce to $\widebar{\vec C}(\pm\matI/2)\vec C^\T = \pm \matI/2$, i.e., the continuity enforcement leaves the second-kind structure of these formulations intact and acts only on the compact contributions.

Notably, for the direct \ac{EFIE} we have two versions corresponding to testing with covariant and contravariant vectors.
Analogously, for the \ac{CFIE} we have included two versions: one discretizing both the \ac{EFIE} and \ac{MFIE} contributions by testing them with contravariant test vectors, and one mixed scheme where the \ac{EFIE} is tested with covariant vectors and the \ac{MFIE} with contravariant vectors. 
In what follows, we refer to the latter scheme as mixed discretization.
It is inspired by the mixed discretization in the context of \ac{MoM} discretizations of the \ac{CFIE} that employ basis and test functions, where \ac{RWG} functions are used to expand the unknown surface density and test the \ac{EFIE}, whereas \ac{BC} functions are used to test the \ac{MFIE}~\cite{begheinSpaceTimeMixedGalerkin2013}.
In our case, the $\veg e_{u/v}$ take the role of primary and the $\veg e^{u/v}$ take the role of dual test functions; however, without incurring any of the overhead commonly associated with dual functions.

    \subsection{Implications}

It should be noted that the proposed explicit continuity-enforcing schemes come with several implications. 
A conforming mesh of the quadrilateral patches is assumed. 
While, in principle, it appears possible to generalize the presented approach to (certain) non-conforming meshes, it requires at least non-trivial bookkeeping, and providing a conforming mesh is a standard feature in many CAD tools and is also common in \ac{MoM} discretizations.
Similarly, as a practically minor drawback, the number of quadrature points in both parametric directions, $u$ and $v$, should be identical. 
For special geometries such as a torus, the proposed scheme allows different numbers of points along toroidal and poloidal lines, but already for a sphere, this is no longer possible. 
At the same time, the gains from the proposed continuity-enforcements outweigh the restriction to conforming meshes and identical quadrature points in $u$ and $v$, as will be shown in the following Section~\ref{resultsSec}.

In terms of computational effort, the three continuity-enforcing matrices are constructed once from the curvilinear geometry description, contain $O(N)$ nonzero entries, and are applied as sparse matrix-vector products before and after each application of the discretized operators.
Their cost is therefore negligible compared to the operator application, and the scheme leaves the latter untouched.

\begin{figure}[tp]
	\centering
	\subfloat[][$\veg j$]{\includegraphics[scale=0.11]{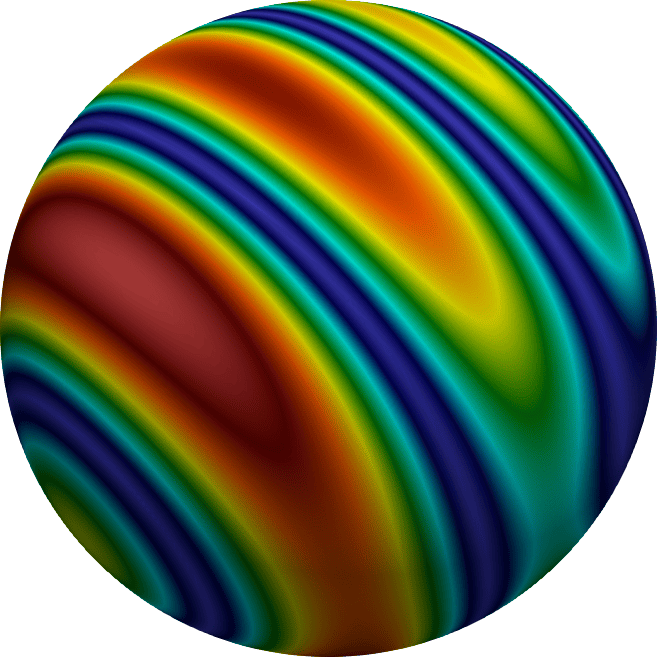}} 	\hspace{1cm}
	\subfloat[][$\veg a$ I-MFIE]{\includegraphics[scale=0.11]{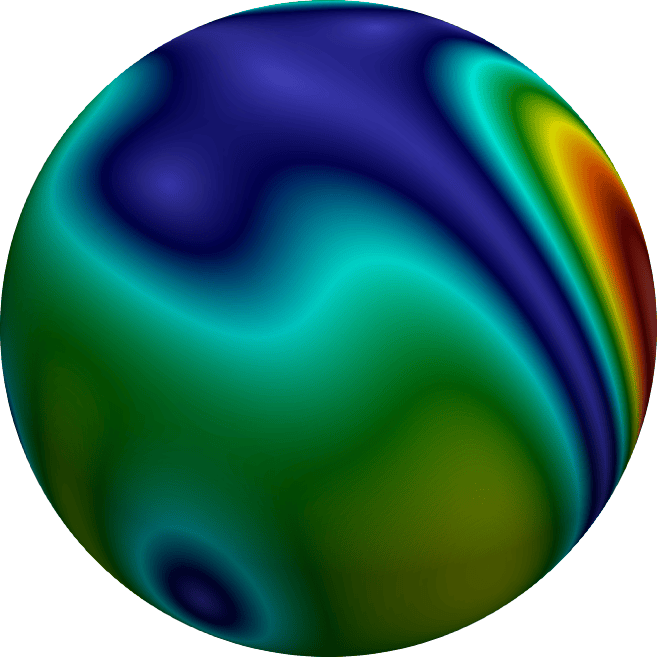}}  \\
	\subfloat[][$\veg b$ I-CFIE]{\includegraphics[scale=0.11]{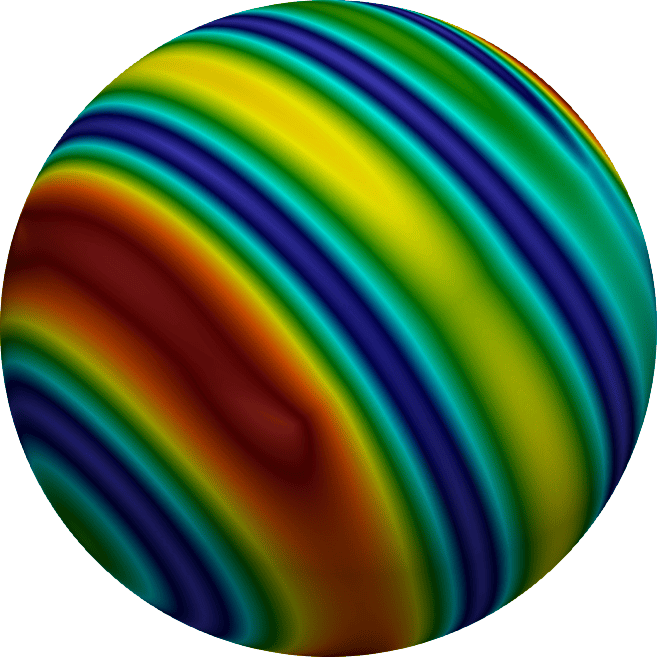}}  \hspace{1cm}
	\subfloat[][$\veg c$ I-CFIE-R]{\includegraphics[scale=0.11]{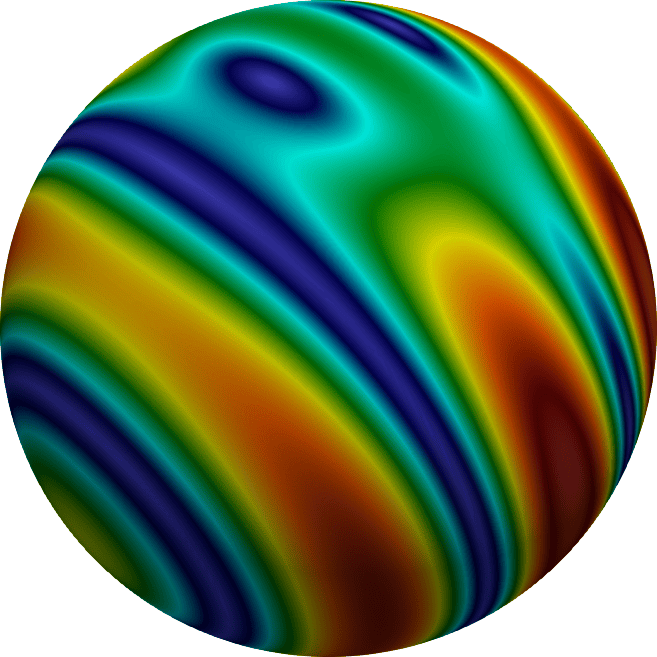}}
	\caption{Induced surface densities $\veg j$, $\veg a$, $\veg b$, and $\veg c$ for a plane wave incident on a sphere of radius $r=\SI{1}{\meter}$ at a frequency of $f=\SI{300}{\mega\hertz}$. The incident angle and the orientation of the spheres are identical.}
	\label{sphereCurrents}
\end{figure} 

\begin{figure*}[tp]
	\centering
	\captionsetup[subfloat]{labelformat=empty}
	\subfloat[][]{\includegraphics[]{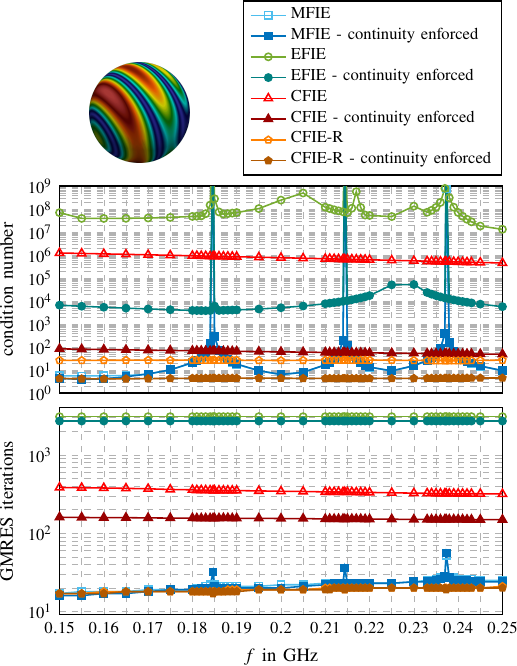}} 	\hfill
	\subfloat[][]{\includegraphics[]{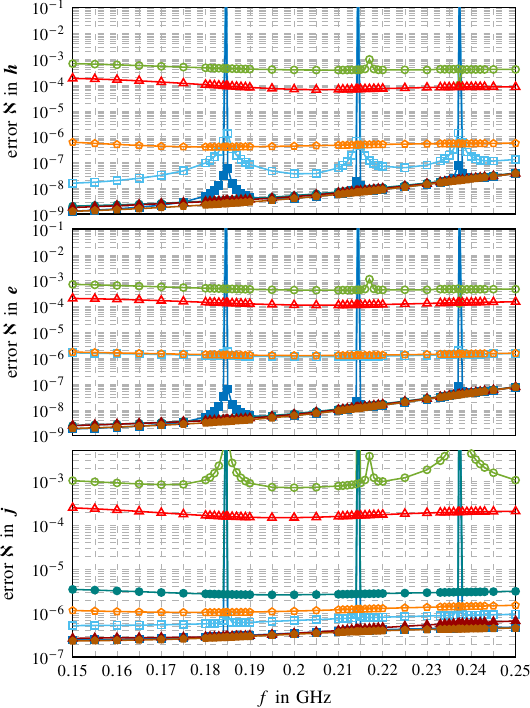}} 
	\caption{Direct formulations: scattering of a plane wave from a sphere of radius $r=\SI{1}{\meter}$ discretized with $N=\num{3072}$ and $N'=\num{2704}$ unknowns without and with continuity enforcement, respectively, for different frequencies. \ac{GMRES} stopping criterion $\epsilon=\num{1e-12}$.}
	\label{sphereCompMie}
\end{figure*} 

\begin{figure*}[tp]
	\centering
	\captionsetup[subfloat]{labelformat=empty}
	\subfloat[][]{\includegraphics[]{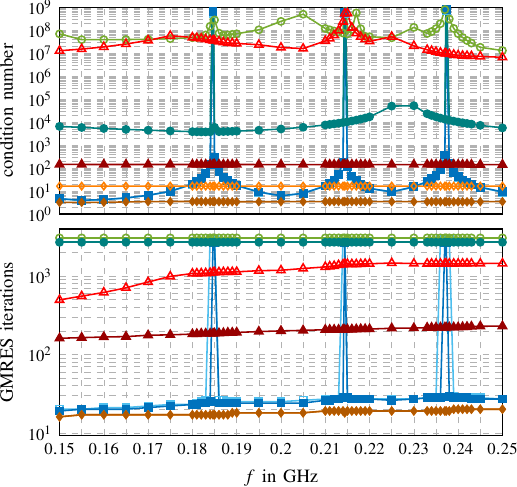}} 	\hfill
	\subfloat[][]{\includegraphics[]{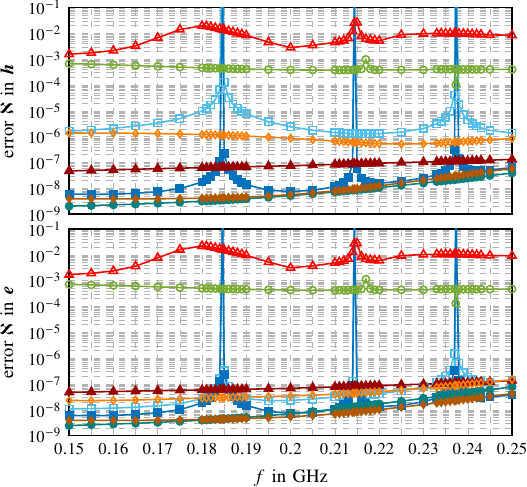}} 
	\caption{Indirect formulations: scattering of a plane wave from a sphere of radius $r=\SI{1}{\meter}$ discretized with $N=\num{3072}$ and $N'=\num{2704}$ unknowns without and with continuity enforcement, respectively, for different frequencies. \ac{GMRES} stopping criterion $\epsilon=\num{1e-12}$.}
	\label{sphereCompIndirect}
\end{figure*} 

\begin{figure*}[tp]
	\centering
	\includegraphics[scale=1]{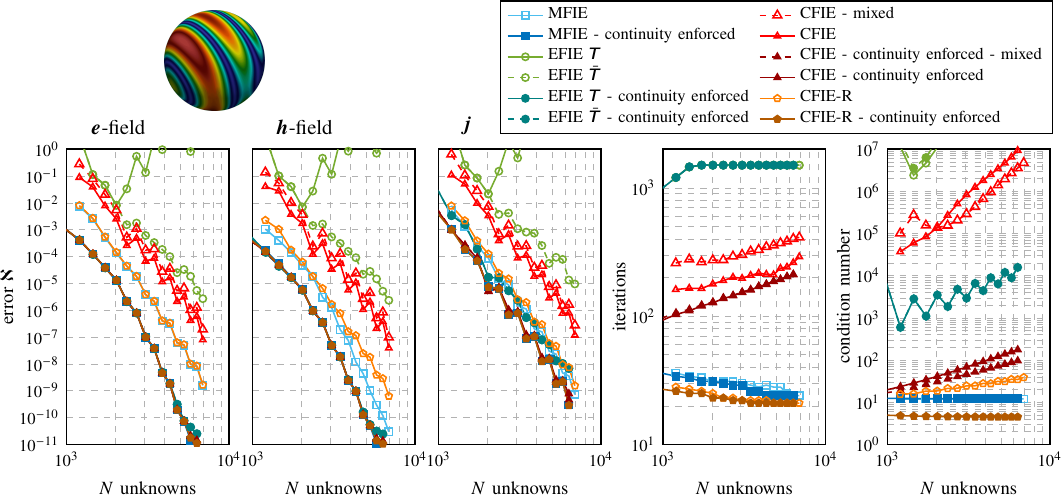}
	\caption{Direct formulations: scattering of a plane wave from a sphere of radius $r=\SI{1}{\meter}$ at a frequency of $f=\SI{300}{\mega\hertz}$ and different numbers of unknowns when refining the number of quadrature points per patch from 10 x 10 to 24 x 24. \ac{GMRES} stopping criterion $\epsilon=\num{1e-12}$.}
	\label{sphereVsN}
\end{figure*}

	\section{Numerical Results} \label{resultsSec}

To investigate the influence of the continuity enforcement on the different integral formulations and the performance of the latter with respect to each other, we consider the scattering of a plane wave from several \ac{PEC} scattering objects: a sphere, a superellipsoid, a spaceplane, and the model of a car.
In all setups, we solve the corresponding \acp{LSE} iteratively employing a \ac{GMRES} solver without restarts.
We study the condition number, the number of iterations, and the achieved error levels.
For the latter, the fields are computed on spherical grids either enclosing or fully inside the objects with $\vartheta$ and $\varphi$ in steps of \SI{5}{\degree}, and the relative average error
\begin{equation}
	\aleph = \sqrt{ \frac{\sum_{ij}{\lVert \veg q(\vartheta_i, \varphi_j) - \hat{\rule{0ex}{0.2ex}\mkern-1.0mu \veg q}(\vartheta_i, \varphi_j)\rVert}_2^2}{\sum_{ij}{\lVert \hat{\rule{0ex}{0.2ex}\mkern-1.0mu \veg q}(\vartheta_i, \varphi_j)\rVert}_2^2} }
	\label{errorDef}
\end{equation}
is computed for $\veg q \in \{\veg e^\mr{sc}, \veg h^\mr{sc}, \veg j \}$ with respect to a reference solution $\hat{\rule{0ex}{0.2ex}\mkern-1.0mu \veg q}$.

    \subsection{Scattering from a Sphere} 

The first scattering object is a sphere of radius $r_\mr{s}=\SI{1}{\meter}$ described exactly (i.e., there is no approximation of the true sphere geometry) by 6 curvilinear patches. 
We employ \mbox{16 x 16} quadrature points per patch, resulting in $N=\num{3072}$ unknowns for the case without continuity enforcement and $N'=\num{2704}$ for the continuity-enforced systems.
The induced surface densities $\veg j$, $\veg a$, $\veg b$, and $\veg c$ are shown in Fig.~\ref{sphereCurrents}.
The different direct formulations are compared in Fig.~\ref{sphereCompMie} in a frequency range of \SI{0.15}{\giga \hertz}-\SI{0.25}{\giga \hertz}
where the reference solutions for the scattered fields and the current density $\veg j$ are computed with respect to a Mie series expansion~\cite{jinTheoryComputationElectromagnetic2015,hofmannSphericalScatteringJuliaPackage2023}) at a distance of $r_\mr{s}=\SI{2}{\meter}$.
The \ac{GMRES} residual is set to $\epsilon = \num{1e-12}$.
For all \ac{EFIE} and \ac{MFIE} formulations, the well-known breakdown (in the condition number, iterations, and errors) can be observed at and close to the theoretical frequencies $f_\mr{R}\approx \SI{0.184662}{\giga\hertz}, \SI{0.214396}{\giga\hertz}, \textrm{and } \SI{0.237299}{\giga\hertz}$ (corresponding to the physical interior resonance frequencies of the \ac{PEC} cavity)~\cite{jinTheoryComputationElectromagnetic2015}.
All \acp{CFIE} overcome the breakdown, removing the peaks in the condition numbers, iterations, and the error in the fields and current densities.

However, there are some more notable differences:
For the \ac{EFIE}, applying the proposed continuity enforcement significantly improves the condition number by four orders of magnitude.
The error in the current density improves by more than two digits, in the electric and magnetic field by about 6 digits, matching the error level of the continuity enforced \ac{MFIE}.
This again shows how crucial continuity enforcement is for the \ac{EFIE}.
The condition number is still notably higher than that of the \ac{MFIE} due to the first kind nature of the integral operator\footnote{Note that applying the regularizer $\vecop R_{\!w}$ to the \ac{EFIE} on its own has no benefit, but, e.g., Calderón preconditioning is required.}.

For the \ac{MFIE}, continuity enforcement does not improve the condition number or the iteration count, but it does improve the error levels:
The error in the surface density $\veg j$ is only slightly improved; however, the error in the magnetic field improves by about a digit, even though the $\vecop T$ operator is not involved at all.
As it is involved in the electric field computation, the error level improves here by up to three digits. 

\begin{figure*}[tp]
	\centering
	\includegraphics[scale=1]{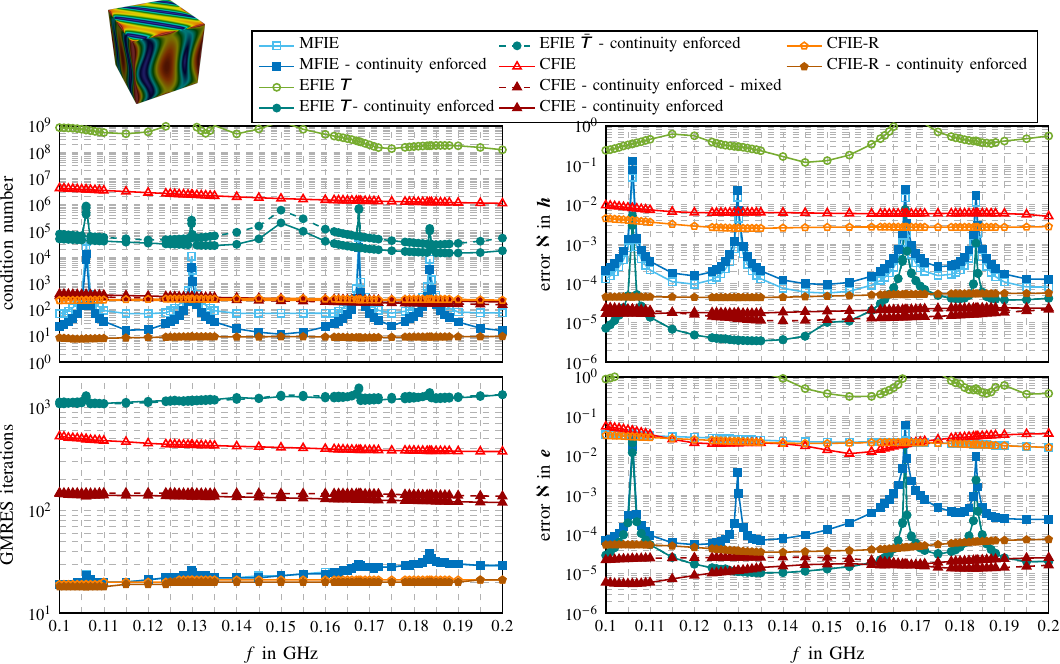}
	\caption{Direct formulations: scattering of a plane wave from a superellipsoid with edge length $\SI{2}{\meter}$ discretized with $N=\num{3888}$ and $N'=\num{3472}$ unknowns without and with continuity enforcement, respectively, for different frequencies. \ac{GMRES} stopping criterion $\epsilon=\num{1e-7}$.}
	\label{cubeEllSineCond}
\end{figure*}
\begin{figure*}[tp]
	\centering
	\includegraphics[scale=1]{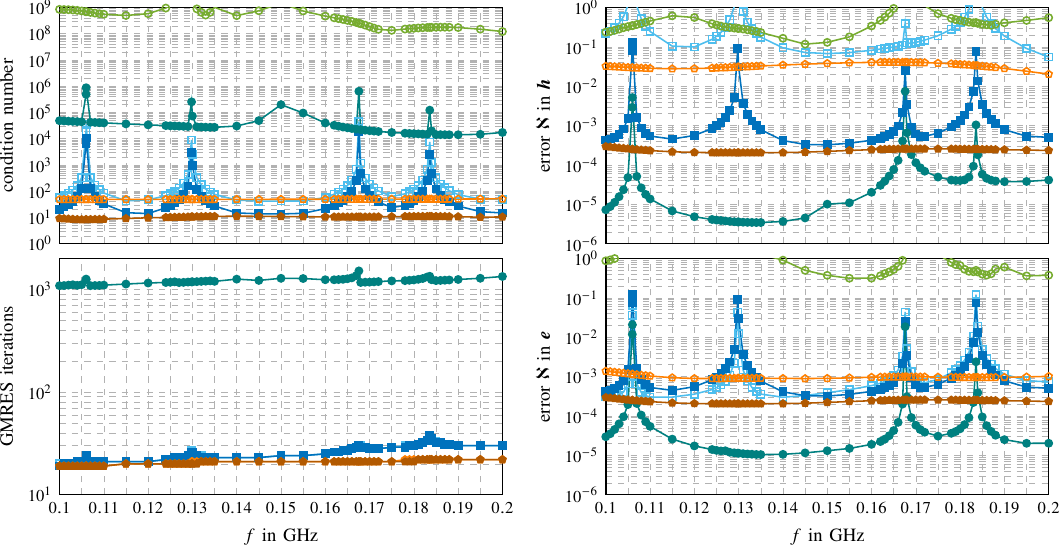}
	\caption{Indirect formulations: scattering of a plane wave from a superellipsoid with edge length $\SI{2}{\meter}$ discretized with $N=\num{3888}$ and $N'=\num{3472}$ unknowns without and with continuity enforcement, respectively, for different frequencies. \ac{GMRES} stopping criterion $\epsilon=\num{1e-7}$.}
	\label{cubeEllSineCondInd}
\end{figure*}
\begin{figure*}[tp]
	\centering
	\includegraphics[scale=1]{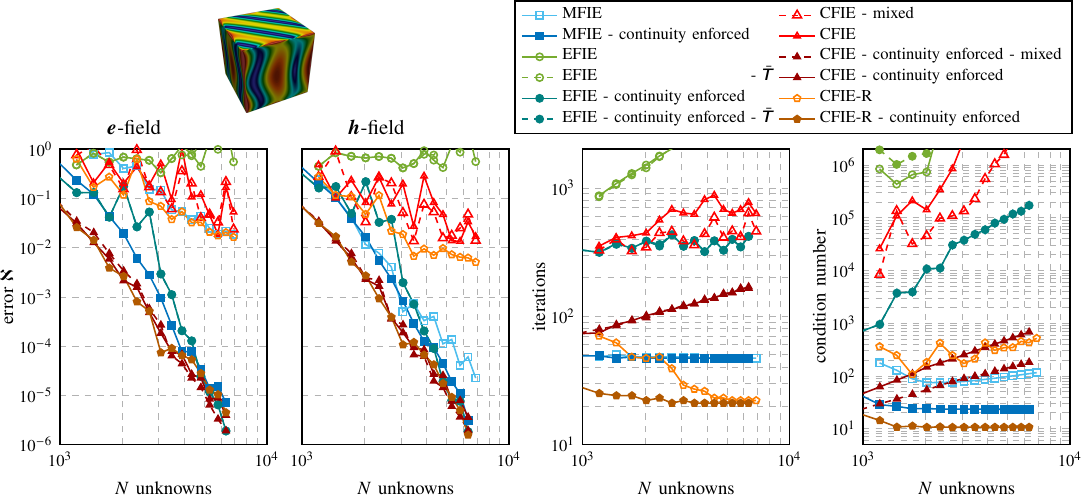}
	\caption{Direct formulations: scattering of a plane wave from a superellipsoid with edge length $\SI{2}{\meter}$ at a frequency of $f=\SI{300}{\mega\hertz}$ and different numbers of unknowns when refining the number of quadrature points per patch from 10 x 10 to 24 x 24. \ac{GMRES} stopping criterion $\epsilon=\num{1e-7}$.}
	\label{superEllVsN}
\end{figure*}

These effects carry over to the \ac{CFIE}, where continuity enforcement improves the condition number, number of iterations, and the error levels significantly, matching the lowest error levels of \ac{EFIE} and \ac{MFIE}, but without breakdowns.
Notably, the regularized \ac{CFIE} already improves the condition number and iteration count even slightly below the continuity-enforced \ac{CFIE}. 
However, the error levels of the latter are not matched, but appear limited by the limited \ac{MFIE} accuracy.
Only the continuity-enforced and regularized \ac{CFIE} achieves the same error level as the continuity-enforced \ac{CFIE} (due to an improvement by another three digits) and the overall lowest condition number and iteration count.

The analogous results for the indirect formulations in Fig.~\ref{sphereCompIndirect} also show the best performance in all aspects\footnote{Note that for the indirect formulations, no reference solution for the surface densities $\veg a$, $\veg b$, and $\veg c$ is available, as the quantities do not represent physical ones.} for the continuity-enforced regularized \ac{CFIE}.
However, there are some notable differences:
Converse to the direct formulation for the \ac{MFIE}, the accuracy improvement by continuity enforcement is stronger in the magnetic than the electric field.
For the \ac{CFIE}, the continuity enforcement appears even more crucial, potentially, due to the $\n \times$ acting on the density before entering the $\vecop T$ operator.
For the regularized \ac{CFIE} it appears to not help that the regularizer $\vecop R_{\!w}$ acts on the unknown density before the surface divergence is computed.

All these effects can also be observed when fixing the frequency to $f=\SI{300}{\mega\hertz}$ for the sphere and studying the convergence when increasing the number of unknowns by increasing the number of quadrature points per patch from 10 x 10 to 24 x 24 as depicted in Fig.~\ref{sphereVsN}.
Here we have also included the case of testing the \ac{EFIE} with contravariant vectors, which for the continuity-enforced variant behaves virtually identically to the \ac{EFIE} tested with covariant vectors.
Moreover, we have included the mixed discretization of the \ac{CFIE}, which shows the same error levels and iteration count as the non-mixed variant. 
Only the condition number seems slightly lower.
This shows that a mixed discretization is possible, but it does not seem to provide a notable advantage.
At the same time, it highlights the flexibility and accuracy of the Nyström approach without requiring a dual basis.

\begin{figure}[tp]
	\centering
	\includegraphics[scale=0.2]{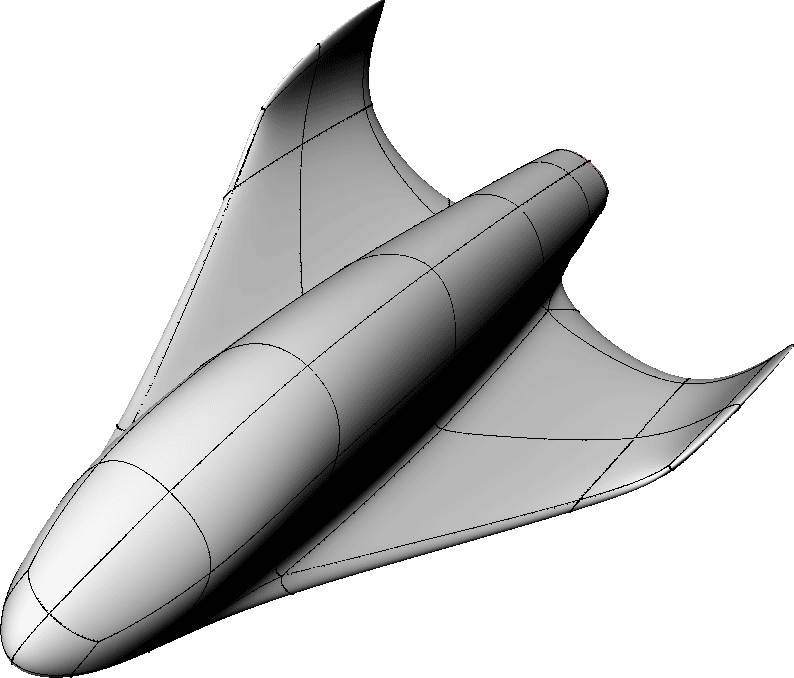}\label{geoS} 	
	\captionsetup[subfloat]{labelformat=empty}
	\caption{Model of a spaceplane of length \SI{13}{\meter} described by \num{92} \ac{NURBS} patches.} 
\label{shuttl}
\end{figure}

\begin{figure*}[tp]
	\centering
	\captionsetup[subfloat]{labelformat=empty}
	\subfloat[][]{\includegraphics[]{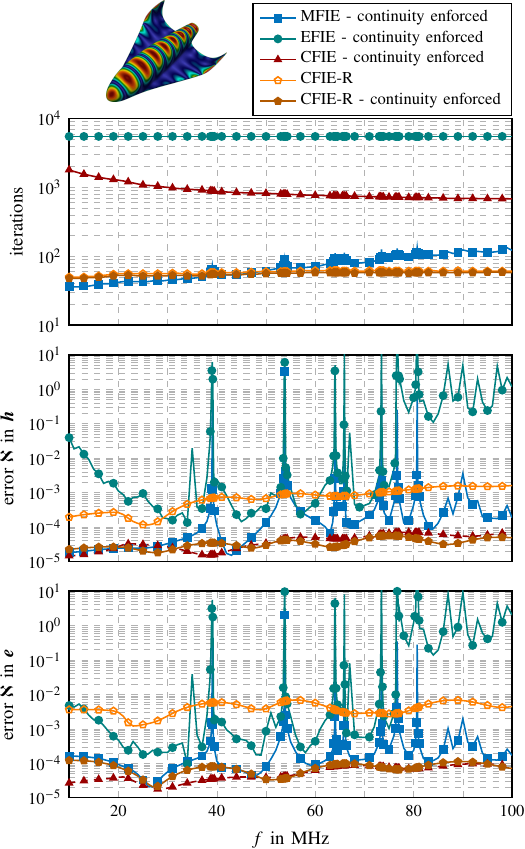}} 	\hfill
	\subfloat[][]{\includegraphics[]{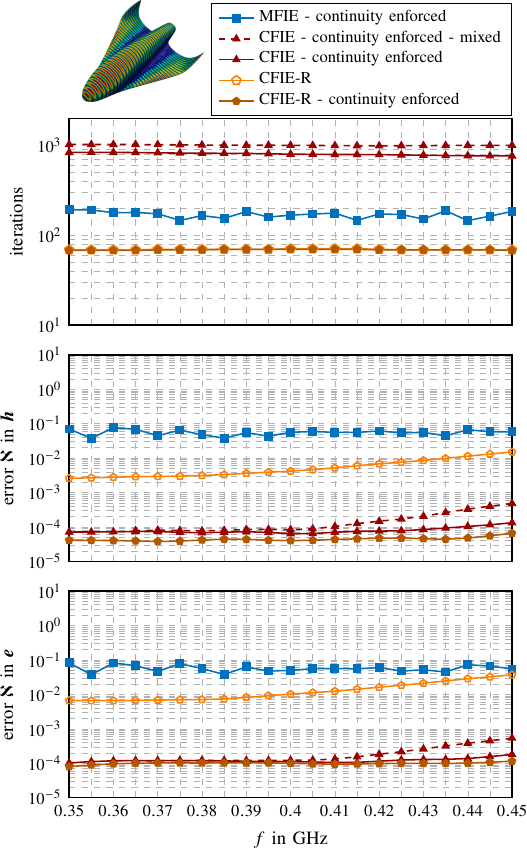}} 
	\caption{Direct formulations: scattering of a plane wave from the spaceplane for different frequencies. Left: discretized with $N=\num{73600}$ and $N'=\num{59620}$ unknowns without and with continuity enforcement, respectively. Right: discretized with $N=\num{16500}$ and $N'=\num{154748}$ unknowns without and with continuity enforcement, respectively. \ac{GMRES} stopping criterion $\epsilon=\num{1e-7}$.}
	\label{shuttle}
\end{figure*} 

\begin{figure*}[tp]
	\centering
	\subfloat[][NURBS model]{\includegraphics[scale=0.17]{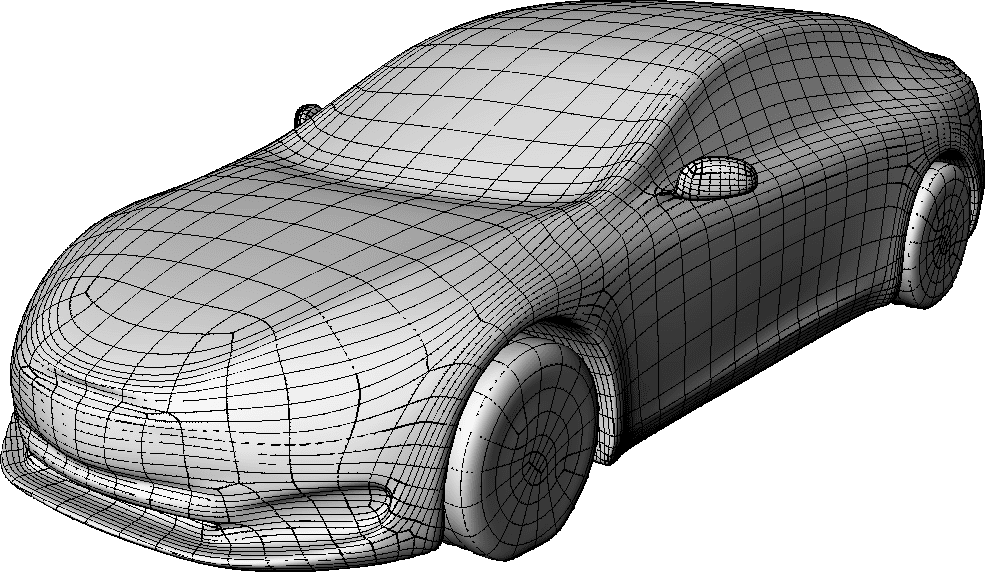}\label{geoCar}} 	\hfill
	\subfloat[][no continuity enforcement]{\includegraphics[scale=0.14]{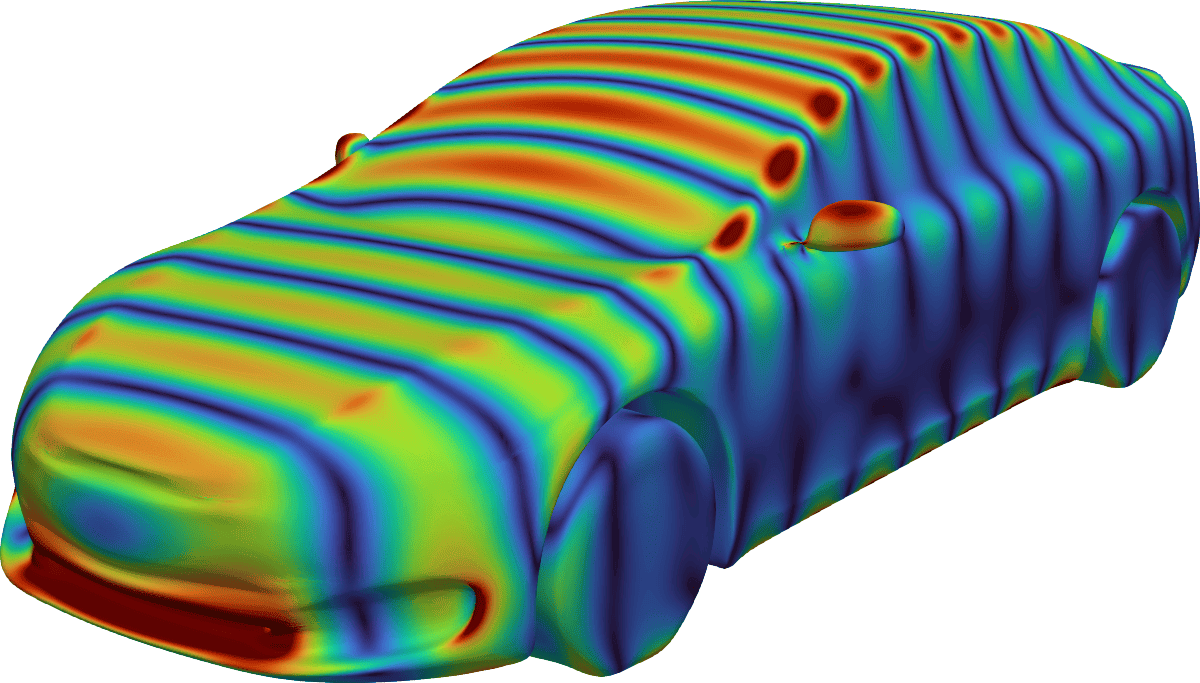}\label{currCar}}  \hfill
	\subfloat[][with continuity enforcement]{\includegraphics[scale=0.14]{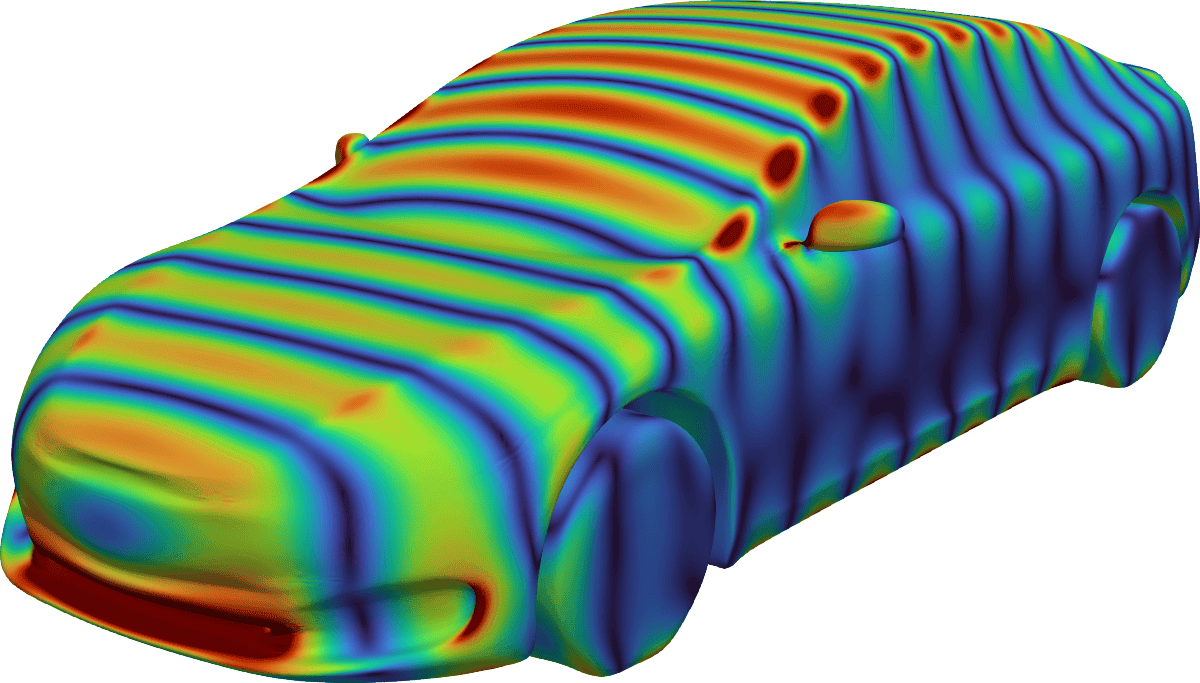}\label{currCarS}}
	\captionsetup[subfloat]{labelformat=empty}
	\caption{Model of a car described by \num{8676} \ac{NURBS} patches and induced surface current densities $\veg j$ for the illumination with a plane wave. The car has a maximum extent of $8.3\lambda$ and is discretized with $N=\num{624672}$ and $N'=\num{433804}$ unknowns without and with continuity enforcement, respectively.}
\label{car}
\end{figure*}
\begin{figure*}[tp]
	\centering
	\subfloat[][no continuity enforcement]{\includegraphics[trim=0 0 0 8, clip, scale=0.16]{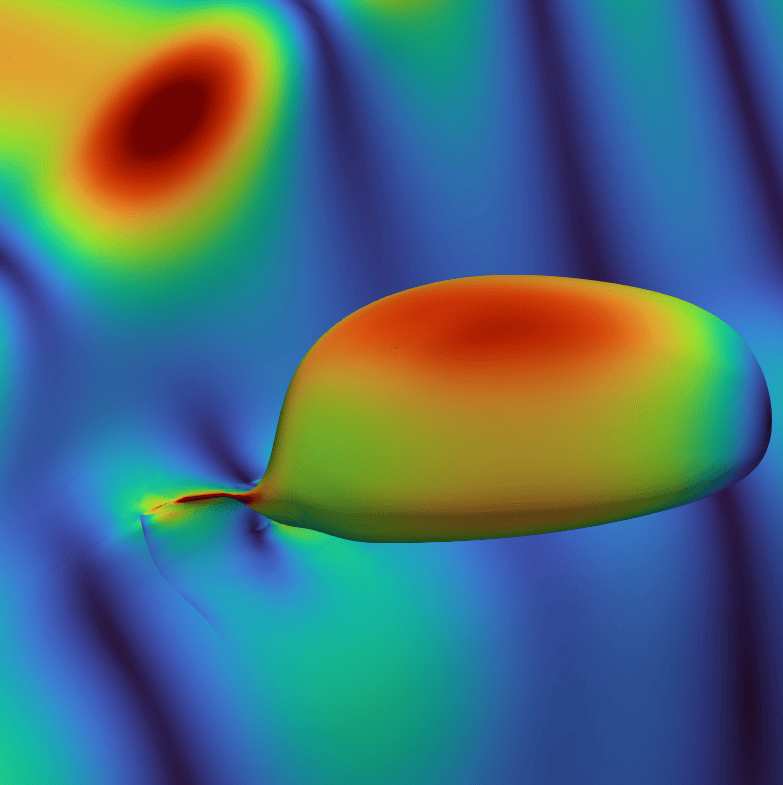}\label{mirr}} 	\hspace{1.5cm}
	\subfloat[][with continuity enforcement]{\includegraphics[trim=0 0 0 0, clip, scale=0.16]{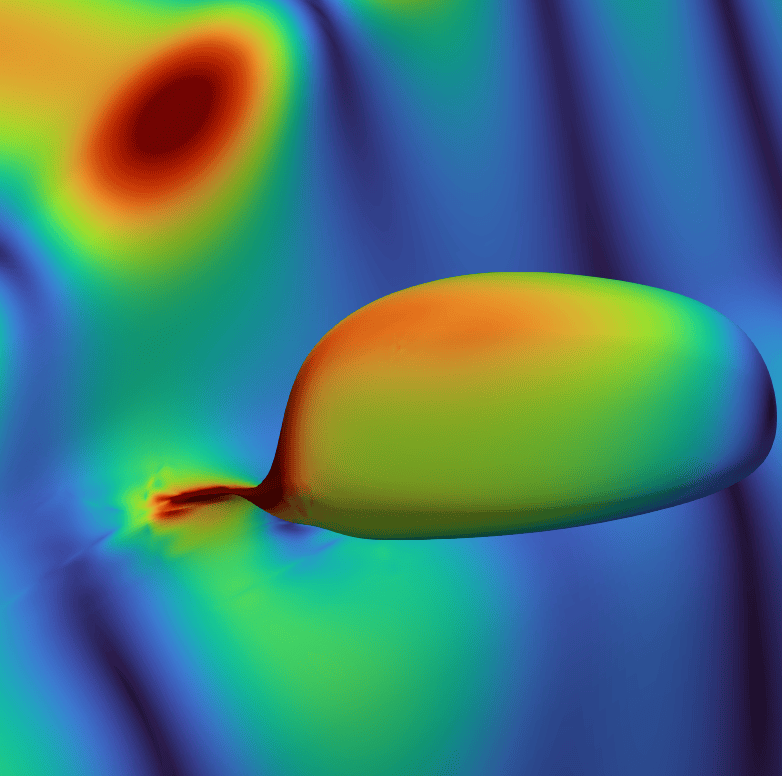}\label{mirrC}}  \hspace{1.5cm}
	\subfloat[][2\textsuperscript{nd}-order MoM solution]{\includegraphics[trim=0 60 0 15, clip, scale=0.16]{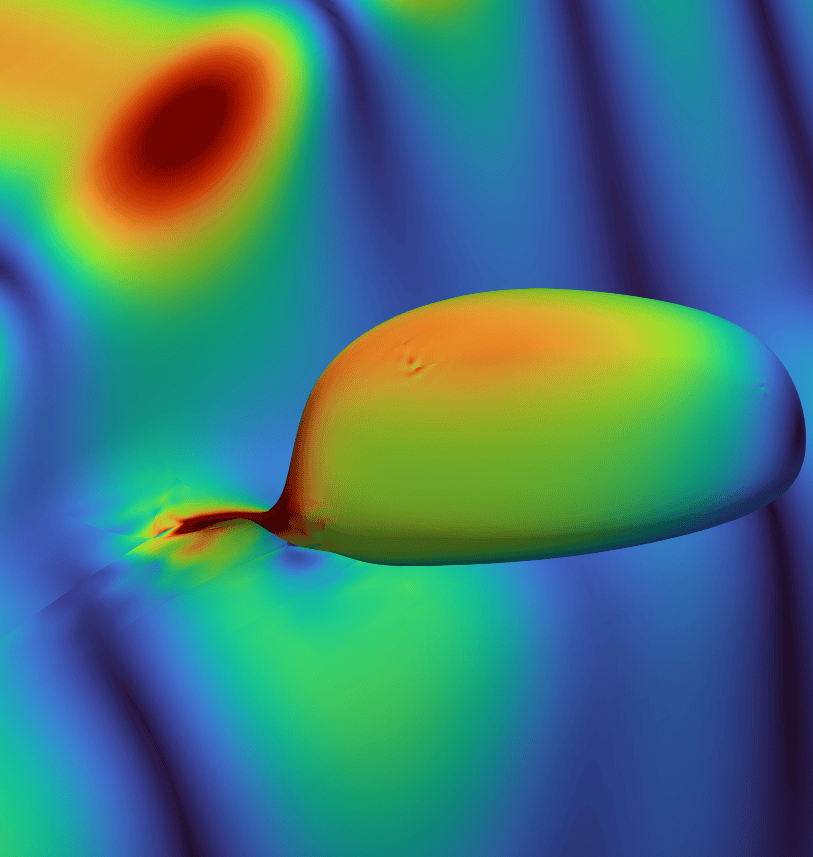}\label{mirrmom}}
	\captionsetup[subfloat]{labelformat=empty}
	\caption{Zoom in on the induced surface current densities on the side mirror: only the continuity enforced density matches the density obtained with a high-order \ac{MoM} solver~\cite{hofmannExplicitHigherOrderDual2026}.}
\label{mirror}
\end{figure*}

    \subsection{Scattering from a Superellipsoid}

To study also the influence of less smooth objects, we consider next a superellipsoid which is implicitly defined by
\begin{equation}
    {\big| x \big|}^n + {\big| y \big|}^n + {\big| z \big|}^n = 1 \,.
\end{equation}
For $n=2$, this corresponds to a sphere, and for $n\rightarrow \infty$, to a cube.
Hence, by setting $n=100$, we obtain an object that is close to a cube, but still smooth. 
We parametrize it via a cube-face projection onto 6 patches, and include a change of variables
\begin{equation}
    u(t) = t + \sin(\uppi t)/(2\uppi)
\end{equation}
to cluster more quadrature points towards the edges.
On this geometry, we employ 18 x 18 quadrature points per patch, resulting in $N=\num{3888}$ unknowns for the case without continuity enforcement and $N'=\num{3472}$ for the continuity-enforced systems.
As no reference solution is available, we determine the error for the direct formulations by computing how well $\veg e^\mr{ex} + \veg e^\mr{sc} = \veg h^\mr{ex} + \veg h^\mr{sc} = \veg 0$ holds on the inside of the superellipsoid.
The \ac{GMRES} residual is set to $\epsilon = \num{1e-7}$.
The results over frequency for the direct formulations are depicted in Fig.~\ref{cubeEllSineCond}.
The overall behavior is similar to the sphere: all \acp{CFIE} overcome the breakdown at resonance frequencies, and continuity enforcement improves conditioning, iterations, and error levels by up to four digits.
Also, in this case, mixed or non-mixed discretization of the \ac{CFIE} makes no notable difference.
However, for this geometry, even the condition number of the \ac{MFIE} is improved by the continuity enforcement.
For the \ac{EFIE}, the latter is even more crucial, translating a virtually non-usable \ac{EFIE} into a very accurate one.
Moreover, the continuity-enforced and regularized \ac{CFIE} is slightly less accurate than the non-regularized continuity-enforced \ac{CFIE}.
Presumably, this is due to the fact that the \ac{EFIE} is more accurate than the \ac{MFIE} in this case, and the regularizer reduces this benefit by converting the \ac{EFIE} contribution to a compact addition to the \ac{MFIE}, whereas the unregularized \ac{CFIE} keeps the unaltered \ac{EFIE}. 

Analogous results can be observed for the indirect formulations in Fig.~\ref{superEllVsN}.
As the total fields are no longer zero in the inside, but exhibit some arbitrary values, we compute a reference solution on a spherical grid of radius $r=\SI{4}{\meter}$ on the outside and compare the fields to the fields obtained by the direct continuity-enforced regularized \ac{CFIE} discretized with 40 x 40 points on each patch corresponding to \num{14704} unknowns.
In addition, in this case, the indirect regularized \ac{CFIE} does not seem to benefit from the regularization being applied to the density before the surface divergence is calculated.

Studying for a fixed frequency of $f=\SI{300}{\mega \hertz}$ the convergence in the number of unknowns in Fig.~\ref{superEllVsN} confirms these results, with the continuity-enforced regularized \ac{CFIE} performing best in all aspects.
Notably, other than \ac{MFIE} and the regularized \ac{CFIE}, only their continuity-enforced counterparts maintain a non-increasing constant condition number, underscoring the benefit of continuity enforcement again.

    \subsection{Scattering from Realistic Geometries}

As a more realistic scatterer, we consider a spaceship model, as shown in Fig.~\ref{shuttl}, described by \num{92} \ac{NURBS} patches. 
We study it in two frequency ranges: $f=\SI{10}{\mega\hertz} - \SI{100}{\mega\hertz}$ and $f=\SI{350}{\mega\hertz} - \SI{450}{\mega\hertz}$.
In the first case, shown in the left of Fig.~\ref{shuttle}, we split each patch into 4 and employ 10 x 10 quadrature points per patch, resulting in $N=\num{73600}$ unknowns for the case without continuity enforcement and $N'=\num{59620}$ for the continuity-enforced systems.
Clearly, the \ac{EFIE} and \ac{MFIE} suffer drastically from resonances, overcome by the \ac{CFIE} and its regularized version.
The latter exhibits a low iteration count (\ac{GMRES} residual $\epsilon=\num{1e-7}$) of 60; and the continuity enforcement improves the error by 2 digits, even though there are \num{13980} fewer unknowns, corresponding to a reduction by approximately 20\,\%.

A similar behavior is observed in the higher frequency range in the right of Fig.~\ref{shuttle}, where 30 x 30 point patches are employed, resulting in $N=\num{165000}$ and $N'=\num{154748}$ unknowns for the cases without and with continuity enforcement, respectively.
In this case, virtually all frequencies are at or near resonance frequencies, and the \ac{MFIE} does not produce accurate results in the whole frequency range.
This is improved by 4 digits of accuracy when using the continuity-enforced regularized \ac{CFIE} despite using \num{10205} fewer unknowns, corresponding to a reduction of about 6\,\%.
Moreover, the low iteration count of 60 is maintained across the entire frequency range.

As a final test case, we consider the CAD model of a car described by \num{8676} \ac{NURBS} patches as shown in Fig.~\ref{car}.
Illuminated at a frequency of $f=\SI{0.5}{\giga\hertz}$, it has a size of about $\SI{5}{\meter} \approx 8.3 \lambda$.
We discretize it with 6 x 6 quadrature points per patch corresponding to $N=\num{624672}$ unknowns for the case without continuity enforcement and $N'=\num{433804}$ for the continuity-enforced systems.
Even though the continuity-enforcement has \num{190868} fewer unknowns, corresponding to a reduction of about 30\,\%, the accuracy improves by 2 digits from \num{1e-1} to \num{1e-3}. 

By closely inspecting the induced surface current densities on the side mirror as depicted in Fig.~\ref{mirror} and comparing them to a high-order \ac{MoM} solution discretized with \num{277632} 2\textsuperscript{nd} order B-spline based basis functions~\cite{hofmannLowFrequencyStabilizationBSpline2024,hofmannExplicitHigherOrderDual2026}, it can be seen that only the continuity-enforced Nyström solution agrees well, further highlighting the advantages of the proposed approach.

    \section{Conclusion}

We have proposed a high-order-accurate continuity-enforcing Nyström-collocation discretization for the direct and indirect \acp{EFIE}, \acp{MFIE}, and (regularized) \acp{CFIE} on smooth surfaces.
The numerical studies of the scattering from a sphere, a superellipsoid, a spaceplane, and a car model confirm that the continuity enforcement proves to be useful beyond discretization of the \ac{EFIE} alone.
For the \ac{EFIE}, it is indispensable: it improves the condition number by up to four orders of magnitude and the error in the scattered fields by up to six digits, turning a virtually unusable formulation into a highly accurate one.
For the \acp{MFIE}, in which the \ac{EFIE} operator does not appear at all, the error levels still improve notably, and the amount of improvement follows the occurrence of that operator in the respective representation formula: it is most substantial for the electric field for the direct and in the magnetic field for the indirect formulations.
Consequently, the common observation that Nyström \acp{MFIE} outperform Nyström \acp{EFIE} holds for the surface density and for the field obtained through the operator $\vecop K$, but not for the field computed through $\vecop T$.

With the accuracy of the \ac{EFIE} restored, \ac{EFIE} and \ac{MFIE} can be combined with equal weights to form \acp{CFIE}, which remove the breakdown at the interior resonance frequencies while inheriting the accuracy of both constituents.
Among all considered formulations, the continuity-enforced regularized \ac{CFIE} performs best, achieving the lowest condition numbers, the lowest iteration counts, and the lowest error levels, all while reducing the total number of unknowns by up to \SI{30}{\percent}.


	
	%

%

	\section*{Acknowledgment}

The authors would like to thank Google Quantum AI for their financial support and Prof. Michel Devoret for many insightful discussions.

	\ifCLASSOPTIONcaptionsoff
	  \newpage
	\fi

	

	\printbibliography

\end{document}